\documentclass[reqno,12pt,a4paper]{amsart}
\usepackage{mathptmx}
\usepackage{graphicx}
\usepackage{amsthm}
\usepackage{amssymb}
\usepackage{color,xcolor}
\usepackage{hyperref}
\usepackage[all,cmtip]{xy}
\usepackage[mathlines]{lineno}

\usepackage[nameinlink]{cleveref}
\crefname{theorem}{Theorem}{Theorems}
\crefname{Maintheorem}{Theorem}{Theorems}
\crefname{lemma}{Lemma}{Lemmas}
\crefname{lem}{Lemma}{Lemmas}
\crefname{remark}{Remark}{Remarks}
\crefname{prop}{Proposition}{Propositions}
\crefname{defn}{Definition}{Definitions}
\crefname{corollary}{Corollary}{Corollaries}
\crefname{cor}{Corollary}{Corollaries}
\crefname{section}{Section}{Sections}
\crefname{subsection}{Subsection}{Sections}
\crefname{figure}{Figure}{Figures}
\crefname{quest}{Question}{Questions}
\crefname{claim}{Claim}{Claims}
\crefname{case}{Case}{Cases}
\crefname{Maincase}{Case}{Cases}
\usepackage{enumitem}
\usepackage{pstricks}
\usepackage{tikz}
 \usepackage{diagbox}
 \usepackage[title]{appendix}
\newtheorem{theorem}{Theorem}[section]
\newtheorem{lemma}[theorem]{Lemma}
\newtheorem{ques}[theorem]{Question}

\newcounter{MainTheoremCounter}

\newtheorem{Maintheorem}[MainTheoremCounter]{Theorem}

\theoremstyle{definition}
\newtheorem{definition}[theorem]{Definition}

\newtheorem{cor}[theorem]{Corollary}
\newtheorem{remark}[theorem]{Remark}

\numberwithin{equation}{section}

\newtheorem{prop}[theorem]{Proposition}

\begin{document}
\title{Saturated theorem along cubes for a measure and applications}

\author[J.~Qiu]{Jiahao Qiu}
\email{qiujh@mail.ustc.edu.cn}

\author[J.~Yu]{Jiaqi Yu}
\email{yjq2020@mail.ustc.edu.cn}

\address[J.~Qiu, J.~Yu]
{Wu Wen-Tsun Key Laboratory of Mathematics, USTC,
Chinese Academy of Sciences and
School of Mathematics, University of Science and Technology of China,
Hefei, Anhui, 230026, P.R. China}

\date{\today}

\begin{abstract}

We show that for a minimal system $(X,T)$, the set of saturated points along cubes with respect to its maximal $\infty$-step pro-nilfactor $X_\infty$
has a full measure. 
As an application, it is shown that
if a minimal system $(X,T)$ has no non-trivial $(k+1)$-tuples with arbitrarily long finite IP-independence sets,
then it has only at most $k$ ergodic measures and
is an almost $k'$ to one extension of $X_\infty$ 
for some $k'\leqslant k$.
Particularly, for $k=1$ we prove that $(X,T)$ is uniquely ergodic (even regular with respect to $X_\infty$), which answers a conjecture stated
in \cite{DDMSY13}.
\end{abstract}


\maketitle

\section{Introduction}
In this section, we will provide the background of the research
and state the main results of the paper.

\subsection{Characteristic factors}\

The study of convergence of the \emph{multiple ergodic averages}
began with Furstenberg's elegant proof of Szemer\'{e}di Theorem \cite{Sz75} via an ergodic theoretical analysis \cite{FH}.
After nearly three decades of work by many researchers,
the problem (for $L^2$ convergence) was finally solved in \cite{HK05,TZ07}.
In their proof, the key idea due to Furstenberg is to find an appropriate factor,
called a \emph{characteristic factor}, that controls the behavior of the limit in the sense of $L^2$ norm.
For the origin of these ideas and this terminology, see \cite{FH}.
The next step is to obtain a specific description for some well-chosen characteristic factors in order to prove convergence.
The result in \cite{HK05,TZ07} shows that
such a characteristic factor can be described as an inverse limit of nilsystems,
which is also called a \emph{pro-nilfactor}.

\medskip

The counterpart of characteristic factors for topological dynamics was first studied by
Glasner in \cite{GE94}. To state the result we need a notion called saturated subset.
Given a map $\pi:X\to Y$ of sets $X$ and $Y$, a subset $L$ of $X$ is called $\pi$-\emph{saturated} if $L=\pi^{-1}(\pi(L))$.
Here is the definition of topological characteristic factors:
  Let $\pi:(X,T)\to (Y,T)$ be a factor map of topological dynamical systems and $d\in \mathbb{N}$.
$(Y,T)$ is said to be a \emph{$d$-step topological characteristic factor}
if there exists a
dense $G_\delta$ subset $\Omega$ of $X$ such that for each $x\in \Omega$ the orbit closure
\[
L_x^d(X):=\overline{\mathcal{O}}\big((\underbrace{x,\ldots,x}_{\text{$d$ times}}),T\times \cdots \times T^d\big)
\]
is $\underbrace{\pi\times \cdots\times \pi}_{\text{$d$ times}}$-saturated.
That is, $(x_1,\ldots,x_d)\in L_x^d(X)$ if and only if
$(x_1',\ldots,x_d')\in L_x^d(X)$ whenever $\pi(x_i)=\pi(x_i')$ for $i=1,\ldots, d$.

In \cite{GE94}, it was shown that for minimal systems, up to a canonically defined proximal
extension, a characteristic family for $T\times \cdots \times T^d$ is the family of canonical PI flows of class $(d-1)$.
In particular, if $(X,T)$ is distal, then its largest class $(d-1)$ distal factor
is its topological characteristic factor along $T\times \cdots \times T^d$. Moreover,
if $(X,T)$ is weakly mixing, then the trivial system is its topological characteristic factor.
For more related results we refer the reader to \cite{GE94}.

\medskip

On the other hand,
to get the corresponding pro-nilfactors for topological dynamics,
in a pioneer work, Host, Kra and Maass \cite{HKM10} introduced the notion of
{\it regionally proximal relation of order $d$}
for a topological dynamical system $(X,T)$, denoted by $\mathbf{RP}^{[d]}$.
For $d\in\mathbb{N}$, we say that a minimal system $(X,T)$ is a \emph{d-step pro-nilsystem}
if $\mathbf{RP}^{[d]}=\Delta$ and this is equivalent for $(X,T)$ being
an inverse limit of minimal $d$-step nilsystems.
For a minimal distal system $(X,T)$, it was proved that
$\mathbf{RP}^{[d]}$ is an equivalence relation and $X/\mathbf{RP}^{[d]}$
is the maximal $d$-step pro-nilfactor \cite{HKM10}.
Later, Shao and Ye \cite{SY12} showed that in fact for
any minimal system the regionally proximal relation of order $d$ is an equivalence
relation.
Furthermore, for any minimal system
$\mathbf{RP}^{[\infty]}=\bigcap_{d\geqslant 1}\mathbf{RP}^{[d]}$
is also an equivalence relation.

\medskip

The recent result in \cite{GHSWY20} improves Glasner's result significantly to pro-nilsystems.
That is, they proved that:

\begin{theorem}[Glasner-Huang-Shao-Weiss-Ye]\label{GHSWY}
Let $(X,T)$ be a minimal system and let $\pi:X\rightarrow X/\mathbf{RP}^{[\infty]}= X_\infty$ be the factor map.
Then there exist minimal systems $X^*$ and $X_\infty^*$ which are almost one to one
extensions of $X$ and $X_\infty$ respectively, and a commuting diagram below such that 
$X_\infty^*$ is a
$d$-step topological characteristic factor of $X^*$ for all $d\geqslant 2$.
\[
\xymatrix{
X \ar[d]_{\pi}  & X^* \ar[d]^{\pi^*}  \ar[l]_{\sigma} \\
X_\infty   & X_\infty^*      \ar[l]_{\tau}
}
\]
\end{theorem}

An equivalent form of \cref{GHSWY} is the following.

\begin{theorem}\label{beta}
Let $(X,T)$ be a minimal system and let $\pi:X\rightarrow X_\infty$
be the factor map.
Then there exists a
dense $G_\delta$ subset $\Omega$ of $X$ such that for every $x\in \Omega$ and every $d\geqslant 2$,
\[
(\pi^{-1}(\pi(x)))^{d}\subseteq  L_x^d(X).
\]
That is, if $\pi(x_i)=\pi(x)$ for $i=1,\ldots, d$, then $(x_1,\ldots,x_d)\in L_x^d(X)$.
\end{theorem}

Motivated by \cref{beta} the following question is natural.

\begin{ques}
Let $(X,T)$ be a minimal system and
let $\pi:X\rightarrow X_\infty$ be the factor map.
Given an ergodic measure $\mu$ on $(X,T)$,
is there a subset $\Omega$ of $X$ with $\mu(\Omega)=1$ such that
 $(\pi^{-1}(\pi(x)))^{d}\subseteq  L_x^d(X)$ for every $x\in \Omega$ and every $d\geqslant 2$?
\end{ques}

We are unable to address this question in the present paper.
Still, we answer a similar question
by considering the face cubegroup instead of the group generated by $T\times T^2\times \cdots \times T^d$,
and derive some interesting applications.

As we know, the face cubegroup introduced in \cite{HKM10}
plays an important role when studying the pro-nilfactors.
See \cref{def-face-group} for the precise definition.
We first study the topological characteristic factors under the action of face cubegroups.
Following a recurrence theorem obtained by Bergelson and McCutcheon \cite{BM00},
and performing a detailed analysis of the Host-Kra cubic measures,
we prove that:

\begin{Maintheorem}\label{full-mesure=1}
Let $(X,T)$ be a minimal system and let $\pi:X\to X/\mathbf{RP}^{[\infty]}$ be the factor map.
Let $\mu$ be an ergodic measure on $(X,T)$ and let
\begin{equation}\label{saturatedpoint}
\Omega=
\{x\in X:
\{x\}\times (\pi^{-1}(\pi(x)))^{2^d-1}\subseteq  \overline{\mathcal{F}^{[d]}}(x^{[d]}),\;\forall d\in\mathbb{N}
\}.
\end{equation}
Then one has $\mu(\Omega)=1$ and
$\Omega$ contains a dense $G_\delta$ set.
\end{Maintheorem}

The point satisfying \eqref{saturatedpoint} will be called a {\it saturated point along cubes}.

\subsection{Applications}\

Following the \emph{local entropy theory}, for a survey see \cite{GY09},
each dynamical system admits a
maximal zero topological entropy factor, and this factor is induced by the smallest
closed invariant equivalence relation containing \emph{entropy pairs} \cite{BL93}.
In \cite{HY06}, entropy pairs are characterized as those pairs that
admit an \emph{interpolating set} of positive density.
Later on, the notions of \emph{sequence entropy pairs} \cite{HLSY03} and \emph{untame pairs}
(called \emph{scrambled pairs} in \cite{WH06}) were introduced.
In \cite{KL07} the concept of \emph{independence}
was extensively studied and used to unify the notions mentioned above.

Let $(X,T)$ be a topological dynamical system and let $\mathcal{U}=(U_1,\ldots,U_k)$ be a tuple
of subsets of $X$. We say that a subset $F\subseteq \mathbb{N}$ is
an \emph{independence set} for $\mathcal{U}$
if for any non-empty finite subset $J\subseteq F$ and any $s=(s(j):j\in J)\in \{1,\ldots,k\}^J$
we have $\bigcap_{j\in J}T^{-j}U_{s(j)}\neq \emptyset$.
It was shown that a tuple of points $x_1,\ldots,x_k$ in $X$ is a sequence entropy $k$-tuple
(if $k=2$ a tuple is said to be a pair) if and only
if each $\mathcal{U}=(U_1,\ldots,U_k)$,
where $U_1,\ldots,U_k$ are neighborhoods of $x_1,\ldots,x_k$ respectively, has arbitrarily long finite
independence sets.
Also, the pair is an untame pair if and only if each $\mathcal{U} = (U_1 ,U_2 )$ as before has infinite independence sets.
It was shown (\cite{EG07,HLSY03,KL07}) that a
minimal null (resp. tame) system is an almost one to one extension of its maximal equicontinuous
factor and is also uniquely ergodic.

For $d\in \mathbb{N}$ and $n_1,\ldots,n_d\in \mathbb{N}$,
we call the set
\[
\mathrm{IP}(n_1,\ldots,n_d):=
\{n_{i_1}+\cdots+n_{i_k}: 1\leqslant i_1<\cdots<i_k\leqslant d,1\leqslant k\leqslant d\}
\]
an \emph{$IP_d$-set}.
The notion of \emph{$Ind_{fip}$-k-tuple} was introduced in \cite{HLY}:
a tuple of points $x_1,\ldots,x_k$ in $X$ is an $\text{Ind}_{fip}$-$k$-tuple if and only
if the independence sets for each $\mathcal{U}=(U_1,\ldots,U_k)$ as before
contain an $\mathrm{IP}_d$-set for every $d\in \mathbb{N}$.
It was showed that
a minimal system without any non-trivial $\mathrm{Ind}_{fip}$-pairs
is an almost one to one extension of
its maximal $\infty$-step pro-nilfactor \cite{DDMSY13}.
They also asked whether a minimal system without any non-trivial $\mathrm{Ind}_{fip}$-pairs is uniquely ergodic?
In this paper, we give a positive answer to this question.

It is also interesting to study the structure of minimal systems without any non-trivial sequence entropy tuples (resp. $\mathrm{Ind}_{fip}$-tuples).
This problem was also studied in \cite{MS07}.
In the recent work, Huang, Lian, Shao and Ye \cite{HLSY21} proved that
if a minimal system has no non-trivial sequence entropy $(k+1)$-tuples
then the system is an almost $k'$ to one extension of its maximal equicontinuous factor for some $k'\in\{1,\ldots,k\}$,
and the system has at most $k^{k}$ ergodic measures.

In this paper, we
study the structure of minimal systems without any non-trivial $\mathrm{Ind}_{fip}$-tuples.
As an application of \cref{full-mesure=1}, we prove that:

\begin{Maintheorem}\label{regular-d}
Let $(X,T)$ be a minimal system without any non-trivial $Ind_{fip}$-$(k+1)$-tuples
and let $\pi:X\to X/\mathbf{RP}^{[\infty]}=Y$ be the factor map.
Then there exists $j\in\{1,\ldots,k\}$ such that $m(Y_{j})=1$, where
$m$ is the unique ergodic measure on $Y$
\footnote{It was proved in \cite{DDMSY13} that $Y$ is an inverse limit of minimal nilsystems, thus $(Y,T)$ is uniquely ergodic.} and 
$Y_i=\{y\in Y: |\pi^{-1}(y)|=i\}, i\in\mathbb{N}$.

Consequently, $X$ has at most $j$ ergodic measures and
$X$ is an almost $j'$ to one extension of $Y$ for some $j'\in\{1,\ldots,j\}$,
that is, $Y_{j'}$ contains a dense $G_\delta$ set.
\end{Maintheorem}

\begin{remark}
Notice that any $\mathrm{Ind}_{fip}$-tuple is a sequence entropy tuple, \cref{regular-d} implies that
if a minimal system has no non-trivial  sequence entropy $(k+1)$-tuples,
then the system has at most $k$ ergodic measures,
which gives a sharp bound on the result obtained in \cite{HLSY21}.
\end{remark}

\begin{remark}
It was shown in  \cite{FGJO21} that every minimal tame system is a regular extension of its maximal equicontinuous
factor  (the set of points with trivial fibres has full measure).
\cref{regular-d} implies that any minimal system without any non-trivial $\mathrm{Ind}_{fip}$-pairs
is a regular extension of its maximal $\infty$-step pro-nilfactor.
\end{remark}


The paper is organized as follows.
In Section 2,
the basic notions used in the paper are introduced.
In Section 3, we study the ergodic decomposition of the Host-Kra cubic measures.
In Section 4, we prove \cref{full-mesure=1}, and as an application we prove  \cref{regular-d} in Section 5.
\medskip

\noindent {\bf Acknowledgments.}
The authors would like to thank Professors Wen Huang, Song Shao and Xiangdong Ye
and Dr. Hui Xu for helping discussions.

\section{Preliminaries}
In this section we gather definitions and preliminary results that
will be necessary later on.
Let $\mathbb{N}$ and $\mathbb{Z}$ be the sets of all positive integers
and integers respectively.
For a finite set $F$, denote by $|F|$ the number of elements of $F$.





\subsection{Dynamical systems}\

A \emph{topological dynamical system} (TDS)  is a pair $(X,T)$,
where $X$ is a compact metric space with a metric $\rho$ and $T : X \to  X$ is a homeomorphism.
For $x\in X,\mathcal{O}(x,T)=\{T^nx: n\in \mathbb{Z}\}$ denotes the \emph{orbit} of $x$.
A TDS $(X,T)$ is called \emph{minimal} if
every point has dense orbit in $X$.
A \emph{homomorphism} between TDSs $(X,T)$ and $(Y,T)$ is a continuous onto map
$\pi:X\to Y$ which intertwines the actions; one says that $(Y,T)$ is a \emph{factor} of $(X,T)$
and that $(X,T)$ is an \emph{extension} of $(Y,T)$. One also refers to $\pi$ as a \emph{factor map} or
an \emph{extension} and one uses the notation $\pi : (X,T) \to (Y,T)$.
An extension $\pi$ is determined
by the corresponding closed invariant equivalence relation
$R_\pi=\{(x,x')\in X\times X\colon \pi(x)=\pi(x')\}$.

\medskip

A {\it measure preserving probability system} (MPS) is a quadruple $(X, \mathcal{X} , \mu, T)$,
where $(X, \mathcal{X} , \mu)$ is a standard Borel probability space
(in particular $X$ is a Polish space and $\mathcal{X}$ is its Borel $\sigma$-algebra) and $T$
is an invertible Borel measure-preserving map ($\mu(TA)=\mu(A)$ for all $A\in \mathcal{X}$ ).
A MPS $(X, \mathcal{X} , \mu,T)$ is {\it ergodic} if for every set $A\in \mathcal{X}$ such
that $TA=A$, one has $\mu(A)=0$ or 1.

\medskip


Let $(X,T)$ be a TDS and $M(X)$ be the set of Borel probability
measures on $X$. Denote by $M(X,T)$ be the set of invariant probability
measures.
For $x_0\in X$ and $\mu\in M(X,T)$, we say that $x_0$ is {\it generic for $\mu$} if
\[
\lim_{N\to \infty}
\frac{1}{N}\sum_{n=1}^{N}f(T^nx_0)=\int_{X}f \mathrm{d}\mu
\]
for any continuous function $f$ on $X$.

\begin{prop}\cite[Proposition 3.7]{FH81}\label{generic-points}
 Let $(X,T)$ be a TDS and $\mu\in M(X,T)$ be ergodic.
 Then almost every point of $X$ with respect to $\mu$ is generic for $\mu$.
\end{prop}

\subsection{Discrete cubes}\

Let $d \geqslant 1$ be an integer and write $[d] = \{1, \ldots , d\}$.
We view the element in $\{0, 1\}^d$ in one of two ways, either as a sequence $\epsilon=(\epsilon_1,\ldots,\epsilon_d)$ of 0's and 1's,
or as a subset of $[d]$. A subset $\epsilon$ corresponds to the sequence $(\epsilon_1,\ldots,\epsilon_d)\in \{0, 1\}^d$
such that $i\in\epsilon$ if and only if $\epsilon_i=1$ for $i\in [d]$. For example,
$\vec{0}=(0,\ldots,0)\in \{0, 1\}^d$ is the same
to $\emptyset\subseteq [d]$.
Let $\{0,1\}^d_*=\{0,1\}^d\backslash\{\vec{0}\}=\{0,1\}^d\backslash\{\emptyset\}$.

If $\vec{n} = (n_1,\ldots, n_d)\in \mathbb{Z}^d$ and $\epsilon\in \{0,1\}^d$, we
define
\begin{equation*}
\vec{n}\cdot \epsilon = \sum_{i=1}^d n_i\cdot\epsilon_i .
\end{equation*}
If we consider $\epsilon$ as $\epsilon\subseteq [d]$, then $\vec{n}\cdot \epsilon =\sum_{i\in\epsilon}n_i$.

For $\alpha=(\alpha_1,\ldots,\alpha_n)\in\{0,1\}^n$ and $\beta=(\beta_1,\ldots,\beta_d)\in\{0,1\}^d$, we define
\begin{equation}\label{sum-vector}
\underline{\alpha\beta}=(\alpha_1,\ldots,\alpha_n,\beta_1,\ldots,\beta_d).
\end{equation}

\medskip

Given a set $X$,
we denote $X^{2^d}$ by $X^{[d]}$. A point $\mathbf{x}\in X^{[d]}$ can be written in one of two equivalent
ways, depending on the context:
\[
\mathbf{x}=(x_\epsilon)_{\epsilon\in\{0,1\}^d}=(x_\epsilon)_{\epsilon\subseteq[d]}.
\]
Hence $x_{\emptyset}=x_{\vec{0}}$ is the first coordinate of $\mathbf{x}$.
Sometimes $\mathbf{x}(\epsilon)$ will be used to denote the $\epsilon$-component of $\mathbf{x}$.
For $x\in X$, we write $x^{[d]}=(x,\ldots,x)\in X^{[d]}$.
The \emph{diagonal} of $X^{[d]}$ is $\Delta^{[d]}=\Delta^{[d]}(X)=\{ x^{[d]}:x\in X\}$.
Usually, when $d=1$, we denote the diagonal by $\Delta_X$ or $\Delta$
instead of $\Delta^{[1]}$.
We can isolate the first coordinate,
writing $X^{[d]}_*=X^{2^d-1}$ and writing $\mathbf{x}\in X^{[d]}$
as $\mathbf{x}=(\mathbf{x}(\vec{0}),\mathbf{x}_*)$,
where $\mathbf{x}_*\in X^{[d]}_*$.

\subsection{Dynamical cubespaces}\
\label{def-face-group}

Let $(X,T)$ be a TDS and $d\in \mathbb{N}$.
We define $\mathbf{Q}^{[d]}(X)$ to be the closure in $X^{[d]}$ of elements of the form
\[
(T^{\vec{n}\cdot\epsilon}x)_{\epsilon\in\{0,1\}^d},
\]
where $\vec{n}\in \mathbb{Z}^d$ and $x\in X$.
When there is no ambiguity, we write $\mathbf{Q}^{[d]}$ instead of $\mathbf{Q}^{[d]}(X)$.

\medskip

{\em Face transformations} are defined inductively as follows: Let
$T^{[1]}_1=\mathrm{id} \times T$. If
$\{T^{[d-1]}_j\}_{j=1}^{d-1}$ is defined already, then set
\begin{align*}
 T^{[d]}_j & =T^{[d-1]}_j\times T^{[d-1]}_j,\; j\in [d-1], \\
 T^{[d]}_d & =\mathrm{id} ^{[d-1]}\times T^{[d-1]}.
\end{align*}

It is easy to see that for $j\in [d]$, the face transformation
$T^{[d]}_j : X^{[d]}\rightarrow X^{[d]}$ can be defined by, for
every ${\bf x} \in X^{[d]}$ and $\epsilon\in \{0,1\}^d $,
\[
T^{[d]}_j{\bf x}=
\left\{
  \begin{array}{ll}
    (T^{[d]}_j{\bf x})(\epsilon)=T\mathbf{x}(\epsilon), & \hbox{$ \epsilon_j=1$;} \\
    (T^{[d]}_j{\bf x})(\epsilon)=\mathbf{x}(\epsilon), & \hbox{$\epsilon_j=0$.}
  \end{array}
\right.
\]

The {\em face cubegroup} of dimension $d$ is the group $\mathcal{F}^{[d]}(X)$ of
transformations of $X^{[d]}$ spanned by the face transformations.
The {\em parallelepiped group} of dimension $d$ is the group
$\mathcal{G}^{[d]}(X)$ spanned by the diagonal transformation and the face
transformations. We often write $\mathcal{F}^{[d]}$ and $\mathcal{G}^{[d]}$ instead of
$\mathcal{F}^{[d]}(X)$ and $\mathcal{G}^{[d]}(X)$ respectively.
For convenience, we denote the orbit closure of $\mathbf{x}\in X^{[d]}$
under $\mathcal{F}^{[d]}$ by $\overline{\mathcal{F}^{[d]}}(\mathbf{x})$,
instead of $\overline{\mathcal{O}(\mathbf{x},\mathcal{F}^{[d]})}$.
For $x\in X$,
let $\mathbf{Q}^{[d]}_x(X)=\mathbf{Q}^{[d]}(X)\cap (\{x\}\times X^{2^d-1})$.

\begin{theorem}\cite{SY12}\label{property}
Let $(X,T)$ be a minimal system and $d\in \mathbb{N}$. Then,
\begin{enumerate}
\item $(\mathbf{Q}^{[d]}(X),\mathcal{G}^{[d]})$ is a minimal system;
\item $(\overline{\mathcal{F}^{[d]}}(x^{[d]}),\mathcal{F}^{[d]})$ is minimal for all $x\in X$;
\item $\overline{\mathcal{F}^{[d]}}(x^{[d]})$ is the unique $\mathcal{F}^{[d]}$-minimal
    subset in $\mathbf{Q}^{[d]}_x(X)$ for all $x\in X$.
\end{enumerate}
\end{theorem}

\subsection{Regionally proximality of higher order}\

\begin{definition}\cite{HKM10}\label{def-rp}
Let $(X,T)$ be a TDS and $d\in \mathbb{N}$.
The \emph{regionally proximal relation of order $d$} is the relation $\textbf{RP}^{[d]}$
defined by: $(x,y)\in\textbf{RP}^{[d]}$ (or $\mathbf{RP}^{[d]}(X)$ in case of ambiguity) if
and only if for every $\delta>0$, there
exist $x',y'\in X$ and $n_1,\ldots,n_d\in \mathbb{N}$ such that:
$\rho(x,x')<\delta,\rho(y,y')<\delta$
and
\[
\rho(  T^{n} x', T^{n}  y')<\delta,
\quad \forall\;n\in\mathrm{IP}(n_1,\ldots,n_d).
\]

We say $(X,T)$ is a \emph{$d$-step pro-nilsystem}
if $\mathbf{RP}^{[d]}$ is trivial.
\end{definition}

\begin{theorem}\cite{SY12}\label{cube-minimal}
Let $(X,T)$ be a minimal system and $d\in \mathbb{N}$.
Then
 \begin{enumerate}
\item $(x,y)\in \mathbf{RP}^{[d]}$ if and only if $(x,y,\ldots,y)=(x,y^{[d+1]}_*)\in\mathbf{Q}^{[d+1]}$
if and only if $(x,y,\ldots,y)=(x,y^{[d+1]}_*)\in \overline{\mathcal{F}^{[d+1]}}(x^{[d+1]})$;
\item $\mathbf{RP}^{[d]}$ is a closed invariant equivalence relation.
\end{enumerate}
\end{theorem}

The regionally proximal relation of order $d$ allows us to construct the maximal $d$-step pro-nilfactor
of a minimal system. That is, any factor of order $d$
factorizes through this system.

\begin{theorem}\label{lift-property}\cite{SY12}
Let $\pi :(X,T)\to (Y,T)$ be the factor map of minimal systems and $d\in \mathbb{N}$. Then,
\begin{enumerate}
\item $(\pi \times \pi) \mathbf{RP}^{[d]}(X)=\mathbf{RP}^{[d]}(Y)$;
\item $(Y,T)$ is a $d$-step pro-nilsystem if and only if $\mathbf{RP}^{[d]}(X)\subseteq R_\pi$.
\end{enumerate}

In particular, the quotient of $(X,T)$ under $\mathbf{RP}^{[d]}(X)$
is the maximal $d$-step pro-nilfactor of $X$.
\end{theorem}

It follows from \cref{cube-minimal} that for any minimal system $(X,T)$,
\[
\mathbf{RP}^{[\infty]}=\bigcap_{d\geqslant 1}\mathbf{RP}^{[d]}
\]
is also a closed invariant equivalence relation.

Now we formulate the definition of $\infty$-step pro-nilsystems.

\begin{definition}\cite{DDMSY13}
A \textbf{minimal system} $(X,T)$ is an \emph{$\infty$-step pro-nilsystem},
if the equivalence relation $\mathbf{RP}^{[\infty]}$ is trivial,
i.e., coincides with the diagonal.
We say
the quotient under $\mathbf{RP}^{[\infty]}$
is the {\it maximal $\infty$-step pro-nilfactor} of $X$.
\end{definition}

\begin{prop}\cite[Proposition 8.1.5]{HSY16}\label{RP-infi}
Let $(X,T)$ be a minimal system and $(x,y)\in X\times X\backslash \Delta$.
Then $(x,y)\in \mathbf{RP}^{[\infty]}$ if and only if
for every neighborhood $V$ of $y$ and any $d\in\mathbb{N}$,
there exists an $\mathrm{IP}_d$-set $I_d$
such that $T^nx\in V$ for all $n\in I_d$.
\end{prop}

\subsection{Nilpotent groups, nilmanifolds and nilsystems}\

Let $G$ be a Lie group. Let $G_1 = G$ and $G_k = [G_{k-1}, G]$ for $k \geqslant 2$,
where $[G, H] = \{[g, h] : g \in G, h \in H\}$ and $[g, h] = g^{-1}h^{-1}gh$.
If there exists some $d \geqslant 1$ such that $G_{d+1} = \{e\}$, $G$ is called a \emph{$d$-step nilpotent Lie group}.
We say that a discrete subgroup $\Gamma$ of a Lie group $G$ is \emph{cocompact} if $G/\Gamma$,
endowed with the quotient topology, is compact.
We say that the quotient $X=G/\Gamma$ is a \emph{$d$-step nilmanifold} if $G$ is a $d$-step nilpotent Lie group and $\Gamma$ is a discrete, cocompact subgroup.
The group $G$ acts on $X$ by left translations and we write this action as $(g,x)\mapsto gx$.
Let $\tau\in G$ and $T$ be the transformation $x\mapsto \tau x$ of $X$.
Then $(X,T)$ is called a \emph{d-step nilsystem}.

We also make use of inverse limits of nilsystems and so we recall the definition of an inverse limit
of systems (restricting ourselves to the case of sequential inverse limits).
If $\{(X_i,T_i)\}_{i\in \mathbb{N}}$ are systems with $\text{diam}(X_i)\leqslant 1$ and $\phi_i:X_{i+1}\to X_i$
are factor maps, the \emph{inverse limit} of the systems is defined to be the
compact subset of $\prod_{i\in \mathbb{N}}X_i$
given by $\{(x_i)_{i\in \mathbb{N}}:\phi_i(x_{i+1})=x_i,i\in \mathbb{N}\}$,
which is denoted by $\lim\limits_{\longleftarrow}\{ X_i\}_{i\in \mathbb{N}}$.
It is a compact metric space endowed with the distance $\rho(x,y)=\sum_{i\in \mathbb{N}}1/ 2^i \rho_i(x_i,y_i)$.
We note that the maps $\{T_i\}$ induce a transformation $T$ on the inverse limit.

The following structure theorems characterize pro-nilsystems.

\begin{theorem}\cite{HKM10}\label{description}
A minimal system is a $d$-step pro-nilsystem if and only if it is an inverse limit of $d$-step minimal  nilsystems.
\end{theorem}

\begin{theorem}\cite{DDMSY13}\label{system-of-order}
A minimal system is an $\infty$-step pro-nilsystem if and only if it is an inverse limit of minimal nilsystems.
\end{theorem}

\subsection{Independence}\

The notion of \emph{independence} was firstly introduced and studied in \cite{KL07,KL16}.
It corresponds to a modification of the notion of \emph{interpolating set}
studied in \cite{GW95,HY06}.

\begin{definition}
Let $(X,T)$ be a TDS. Given a tuple $\mathcal{U} = (U_1,\ldots,U_k )$ of subsets
of $X$ we say that a subset $F\subseteq \mathbb{N}$  is an \emph{independence set} for $\mathcal{U}$ if for any non-empty
finite subset $ J\subseteq  F$ and any $s=(s(j):j\in J)\in \{1,\ldots,k\}^J$  we have
\[
\bigcap_{j\in J}T^{-j}U_{s(j)}\neq \emptyset.
\]

We shall denote the collection of all independence sets for $\mathcal{U}$ by $\mathrm{Ind}(U_1 ,\ldots,U_k )$ or $\mathrm{Ind}(\mathcal{U})$.
\end{definition}

The following definition of $\mathrm{Ind}_{fip}$-$k$-tuple comes from \cite{DDMSY13}.

\begin{definition}\label{def-IN-d}
Let $(X,T)$ be a TDS and $k\geqslant 2$.
A tuple $(x_1,\ldots,x_k ) $ of points in $X$ is called an
{\it $Ind_{fip}$-$k$-tuple}
if for any neighborhoods  $U_1,\ldots,U_k$ of $x_1,\ldots,x_k$ respectively and any $d\in \mathbb{N}$,
 there exists an $\mathrm{IP}_d$-set $I_d$ such that
$I_d$
is an independence set for
$(U_1,\ldots,U_k)$.

Denote by $\mathrm{Ind}_{fip}^{(k)} (X)$ the set of all
$\mathrm{Ind}_{fip}$-$k$-tuples of $(X,T)$.
\end{definition}

The following characterization of  $\mathrm{Ind}_{fip}$-tuples using
dynamical cubespaces was first given in \cite{DDMSY13} for $k=2$.

\begin{lemma}\label{tuple}
Let $(X,T)$ be a minimal system, $k\geqslant 2$ and $x_1,\ldots,x_k\in X$. Then $\{x_1,\ldots,x_k\}^{[d]}\subseteq\mathbf{Q}^{[d]}(X)$ for all $d\in\mathbb{N}$ if and only if $(x_1,\ldots,x_k)$ is an $Ind_{fip}$-$k$-tuple.
\end{lemma}

\begin{proof}
Let $x_1,\ldots,x_k\in X$ be distinct points.
First we assume $\{x_1,\ldots,x_k\}^{[d]}\subseteq \mathbf{Q}^{[d]}(X)$
for all $d\in \mathbb{N}$.

Fix $d\in\mathbb{N}$.
Let $ \{1,\ldots,k\}^{[2^d]} =\{\omega_1,\ldots,\omega_{k^{2^d}}\}$.
For $1\leqslant i\leqslant k^{2^d}$, set $\omega_i=(\omega_{i,\epsilon})_{\epsilon\in \{0,1\}^d}$.
Choose $m\in\mathbb{N}$ with $2^m\geqslant k^{2^d}$ and
let $\omega_i=\vec{1}\in\{0,1\}^{2^d}$ for $k^{2^d}+1\leqslant i\leqslant 2^m$.
Let $\{0,1\}^m=\{\alpha_1,\ldots,\alpha_{2^m}\}$.
Notice that for every $\eta\in \{0,1\}^{m+d}$, there exists unique $i\in\{1,\ldots,2^m\}$ and $\beta\in \{0,1\}^d$ such that
$\eta=\underline{\alpha_i\beta}$\footnote{See \eqref{sum-vector} for the definition of $\underline{\alpha_i\beta}$.}.

Let $U_1,\ldots, U_k$ be neighborhoods of $x_1,\ldots, x_k$ respectively.
Now we define $\mathbf{x}=(x_\eta)_{\eta\in \{0,1\}^{m+d}}\in \{x_1,\ldots,x_k\}^{2^{m+d}}$ and
$\mathbf{U}=(U_\eta)_{\eta\in \{0,1\}^{m+d}}\in \{U_1,\ldots,U_k\}^{2^{m+d}}$
such that
$x_\eta=x_{\omega_{i,\beta}}$ and $U_\eta=U_{\omega_{i,\beta}}$ if $\eta=\underline{\alpha_i\beta}$.
Then $\mathbf{U}$ is a neighborhood of $\mathbf{x}$ and $\mathbf{x}\in \mathbf{Q}^{[m+d]}(X)$
by our assumption.
Thus by \cref{property} there exist $a\in\mathbb{N}$ and $\vec{a}=(a_1,\ldots,a_m,n_1,\ldots,n_{d})\in\mathbb{N}^{m+d}$ such that
$T^{a+\vec{a}\cdot\eta}x_1\in U_{\eta}$ for all $\eta\in \{0,1\}^{m+d}$.

We claim that the $\mathrm{IP}_d$-set generated by $n_{1}, \ldots , n_{d}$ belongs to $\mathrm{Ind}(U_1,\ldots,U_k)$.
Aa a matter of fact, for any $s: \{0,1\}^d\to \{1,\ldots,k\}$
there exists $\omega_i$ such that $\omega_{i,\epsilon}=s(\epsilon)$ for all $\epsilon\in \{0,1\}^d$.
Set $\vec{b}=(a_1,\ldots,a_m)$ and $\vec{n}=(n_1,\ldots,n_d)$.
Then for any
$\epsilon\in \{0,1\}^d$ we have
\[
T^{\vec{n}\cdot\epsilon}(T^{a+\vec{b}\cdot \alpha_i}x_1)=
T^{a+\vec{a}\cdot \underline{\alpha_i\epsilon}}x_1\in U_{\underline{\alpha_i\epsilon}}=
U_{\omega_{i,\epsilon}}=U_{s(\epsilon)}
\]
which implies that
\[
T^{a+\vec{b}\cdot \alpha_i}x_1\in
\bigcap_{\epsilon\in \{0,1\}^d} T^{-\vec{n}\cdot\epsilon}U_{s(\epsilon)}.
\]
As $d$ is arbitrary,
we conclude that  $(x_1,\ldots,x_k)$ is an $\mathrm{Ind}_{fip}$-$k$-tuple.

\medskip

Conversely. Assume that $(x_1,\ldots,x_k)$ is an $\mathrm{Ind}_{fip}$-$k$-tuple.
That is, for any neighborhoods $U_1,\ldots, U_k$ of $x_1,\ldots, x_k$ respectively and any $d\in\mathbb{N}$, there exists
$\vec{m}_0=(m_1,\ldots,m_d,m_{d+1})\in \mathbb{N}^{d+1}$ such that for any $t\in\{1,\ldots,k\}^{\{0,1\}^{d+1}}$,
\[
\bigcap_{\omega\in \{0,1\}^{d+1}}T^{-\vec{m}_0\cdot \omega}U_{t(\omega)}\neq\emptyset.
\]

Now fix $s\in\{1,\ldots,k\}^{\{0,1\}^d}$.
There exists  $t\in\{1,\ldots,k\}^{\{0,1\}^{d+1}}$ such that $t(\underline{\epsilon1})=s(\epsilon)$ for all $\epsilon\in \{0,1\}^d$.
Let $x\in \bigcap_{\omega\in \{0,1\}^{d+1}}  T^{-\vec{m}_0\cdot \omega}U_{t(\omega)}$ and set $\vec{m}=(m_1,\ldots,m_d)$.
Then we have
\[
T^{\vec{m}\cdot\epsilon}(T^{m_{d+1}}x)=T^{\vec{m}_0\cdot \underline{\epsilon1}}
\in U_{t(\underline{\epsilon1})}=U_{s(\epsilon)}
\] for all $\epsilon\in \{0,1\}^d$ which implies
 $\{x_1,\ldots,x_k\}^{[d]}\subseteq \mathbf{Q}^{[d]}(X)$.

 This completes the proof.
\end{proof}

\subsection{Some facts about hyperspaces.}\

Let $X$ be a compact metric space with a metric $\rho$.
Let $2^X$ be the collection of non-empty closed subsets of $X$.
We may define a metric on $2^X$ as follows:
\begin{align*}
\rho_H(A,C) &= \inf\{  \epsilon>0:A\subseteq B(C,\epsilon),C\subseteq B(A,\epsilon) \} \\
& =\max\{ \max_{a\in A} \rho(a,C),\max_{c\in C} \rho(c,A)\},
\end{align*}
where $\rho(x,A)=\inf_{y\in A}\rho(x,y)$ and $B(A,\epsilon)=\{x\in X:\rho(x,A)<\epsilon\}$.
The metric $\rho_H$ is called the \emph{Hausdorff metric} on $2^X$.

\medskip

Let $X,Y$ be two compact metric spaces.
Let $F : Y \to 2^X$ be a map and $y\in Y$.
We say that $F$ is {\it lower semi-continuous} ({\it l.s.c.}) at $y$ if for any $\epsilon>0$ there exists a neighbourhood $U$ of $y$
such that $F(y) \subseteq B(F(y'),\epsilon)$ for all $y'\in U$;
and $F$ is {\it upper semi-continuous} ({\it u.s.c.}) at $y$
 if for any $\epsilon>0$ there exists a neighbourhood $U$ of $y$
such that $F(U) \subseteq B(F(y),\epsilon)$.

We have the following well known result, for a proof see
\cite[Theorem 2.13]{KK66} and \cite[p.70-71]{KK68}.

\begin{theorem}\label{con-usc}
 Let $X,Y$ be compact metric spaces. If $F : Y \to 2^X$ is u.s.c., then
the points of continuity of $F$ form a dense $G_\delta$ set in $Y$.
\end{theorem}

\subsection{Fundamental extensions}\

Let $\pi:(X,T)\to (Y,T)$ be a factor map of TDSs.
One says that:
\begin{enumerate}
  \item  $\pi$ is an {\it open} extension if it is open as a map;
  \item  $\pi$ is an {\it almost $k$ to one} extension if
  there exists a dense $G_\delta$ subset $\Omega$ of $Y$
  such that $|\pi^{-1}(y)|=k$ for every $y\in \Omega$.
\end{enumerate}

The following is a well known fact about open mappings (see \cite[Appendix A.8]{JDV} for example).

\begin{theorem}\label{open-map}
Let $\pi:(X,T)\to (Y,T)$ be a factor map of TDSs.
Then the map $\pi^{-1}:Y\to 2^X,y \mapsto  \pi^{-1}(y)$ is continuous
if and only if $\pi$ is open.
\end{theorem}

\begin{prop}\cite[Proposition 2.14]{HLSY21}\label{almost-finite-to-one}
Let $\pi:(X,T)\to (Y,T)$ be a factor map of TDSs with $(Y,T)$ being minimal.
If $Y_f=\{y\in Y:|\pi^{-1}(y)|<\infty\}$ is non-empty,
then $Y_0=\{y\in Y:|\pi^{-1}(y)|=N\}$ contains a dense $G_\delta$ subset of $Y$,
where $N=\min_{y\in Y_f}|\pi^{-1}(y)|$.
\end{prop}

\subsection{The Host-Kra cubic measures}\

We recall the definition of the dual function from \cite{HK18,HKM10}.

\begin{definition}

Let $(X,\mathcal{X},\mu,T)$ be an ergodic system and $d\in\mathbb{N}$.
For every $f\in L^\infty(\mu)$, the limit function
\[
\lim_{N\to \infty}\frac{1}{N^d}\sum_{\vec{n}\in [N]^d}\prod_{\epsilon\in \{0,1\}^d_*}f(T^{\vec{n}\cdot \epsilon}x)
\]
is called the {\it dual function of order d of f} and is written $\mathcal{D}_df$.
  
\end{definition}

\begin{lemma}\cite[Lemma 5.4]{HKM10}\label{dual-function}
If $(X,\mathcal{X},\mu,T)$ is an ergodic system and $d\geqslant 1$ is an integer,
then for every $A\in \mathcal{X}$ we have $\mathcal{D}_d\mathbf{1}_A(x)>0$ for $\mu$-a.e. $x\in A$.
\end{lemma}



For $d\in\mathbb{N}$  and $k\in[d]$,
we say a map $p:\{0,1\}^k\to \{0,1\}^d$ 
is a {\it morphism} if there exist  pairwise disjoint non-empty subsets
$\epsilon_1,\ldots,\epsilon_k$ of $[d]$ such that $p(\emptyset)=\emptyset$
and  $p(\eta)=\bigcup_{i\in\eta}\epsilon_i$ for $\eta\in\{0,1\}^k_*$.
Then $p$ determines a projection $\bar{p}:X^{[d]}\to X^{[k]}$ such that
\[
\bar{p}(\mathbf{x})(\eta)=\mathbf{x}(p(\eta))
\]
for $\mathbf{x}\in X^{[d]}, \eta\in\{0,1\}^k$.

\medskip

The following theorem comes from  \cite{HSY17}.

\begin{theorem}\label{supp-order-d}
Let $(X,T)$ be a minimal system and $\mu\in M(X,T)$ be an ergodic measure.
Then there exists a family $\{\mu_x^{[d]}=\delta_x\times \mu_x^{[d]_*} \}_{x\in X,d\in\mathbb{N}}$ of Borel probability measures on $X^{[d]}$
with $\mu_x^{[1]_*}=\mu$ such that
the following statements hold.
\begin{enumerate}
\item  $\mu_x^{[d]}$ is $\mathcal{F}^{[d]}$-ergodic for $\mu$-a.e. $x\in X$.
\item For any $k\in[d]$ and any morphism $p:\{0,1\}^k\to \{0,1\}^d$, one has $\bar{p}_*(\mu_x^{[d]})=\mu_x^{[k]}$
     for $\mu$-a.e. $x\in X$.
\item 
For any $f\in L^\infty(\mu)$ one has
\[
\mathcal{D}_df(x)=\int_{X_*^{[d]}}\prod_{\epsilon\in \{0,1\}^d_*}f\mathrm{d}\mu_x^{[d]_*}.
\]

  \end{enumerate}
\end{theorem}

\begin{prop}\label{support-d}
Let $(X,T)$ be a minimal system and $\mu\in M(X,T)$ be an ergodic measure.
Let $\{\mu_x^{[d]}\}_{x\in X,d\in\mathbb{N}}$ be as in \cref{supp-order-d}
and let $M$ be a measurable subset of $X$ with positive measure.
Then there exists a measurable subset $M_0$ of $M$ with $\mu(M_0)=\mu(M)$ such that
for any $d\in\mathbb{N},x\in M_0$ and any neighborhood $V$ of $x$, one has
\[
\mu_x^{[d]}((M\cap V)^{[d]})>0.
\]
\end{prop}

\begin{proof}
Let $\{V_i\}_{i\in\mathbb{N}}$ be a countable base for the topology of $X$. 
For every $i\in\mathbb{N}$,
it follows from \cref{dual-function} that there exists a measurable subset $V_i'\subseteq M\cap V_i$ with $\mu(V_i')=\mu( M\cap V_i)$ such that
$\mathcal{D}_d\mathbf{1}_{M\cap V_i}(x)>0$ for $x\in V_i'$.
Thus by \cref{supp-order-d}  for any $x\in V_i'$ we have
\[
\mu_x^{[d]} ((M\cap V_i)^{[d]})=\delta_x(M\cap V_i)\times \mu_x^{[d]_*}((M\cap V_i)^{2^d-1})
=\mathcal{D}_d\mathbf{1}_{M\cap V_i}(x)>0.
\]
By setting $M_d=\bigcap_{i\in\mathbb{N}}V_i'\cup (M\backslash V_i )$
and $M_0=\bigcap_{d\in\mathbb{N}}M_d$ we have $\mu(M_0)=\mu(M)$.
We next show that $M_0$ meets our requirement.

Let $d\in\mathbb{N},x\in M_0$ and $V$ be a neighborhood of $x$.
Recall that  $\{V_i\}_{i\in\mathbb{N}}$ is a base for the topology of $X$,
we can choose $V_i$ with $x\in V_i\subseteq V$.
This implies that
\[
\mu_x^{[d]} ((M\cap V)^{[d]})\geqslant
\mu_x^{[d]} ((M\cap V_i)^{[d]})
>0.
\]

This completes the proof.
\end{proof}

\section{Proof of \cref{full-mesure=1}}

In this section, we prove \cref{full-mesure=1}.
Throughout this section, we assume that $(X,T)$ is a minimal system and $\mu\in M(X,T)$ is an ergodic measure.
Let $\pi:X\to Y=X/\mathbf{RP}^{[\infty]}$ be the factor map.


\subsection{An intermediate proposition of proving \cref{full-mesure=1}}\

We start with
the following definition.

\medskip

\noindent {\bf Definition.}
Let $M$ be a closed subset of $X$ with positive measure.
We say that
the point $(x_\epsilon)_{\epsilon\in \{0,1\}^d}\in X^{[d]}$ has \textbf{the property $S_d(M)$}
 if
$x_{\vec{0}}=x\in M$, $\pi(x_\epsilon)=\pi(x),\epsilon\in \{0,1\}^d$ and
for any $\delta>0$
there exists a closed set $W\subseteq M\cap B(x,\delta)$
with positive measure and $\vec{n}\in\mathbb{Z}^d$ such that
\[
W\subseteq \bigcap_{\epsilon\in \{0,1\}^d}T^{-\vec{n}\cdot \epsilon} B_M(x_\epsilon,\delta),
\]
where
\[
B_M(x_\epsilon,\delta)=
\begin{cases}
M\cap B(x,\delta), &\mathrm{if}\;x_\epsilon=x, \\
B(x_\epsilon,\delta),& \mathrm{otherwise}.
\end{cases}
\]

For $x\in X$ and $d\in \mathbb{N}$, let
\[
A_x^{[d]} = \{x\}\times (\pi^{-1}(\pi(x)))^{2^d-1}.
\]

\bigskip

To show \cref{full-mesure=1}, it suffices to show \cref{order-d}.
As the proof of this implication is clear, we will give the proof of  \cref{full-mesure=1} assuming \cref{order-d} in  the remainder of this subsection.

We remark that the proof of \cref{order-d} is very long and we will divide the proof into three parts:
First we show that there exists $x\in M$ such that if $(x,y)\in \mathbf{RP}^{[\infty]}$,
then $(x,y,x^{(2^d-2)})$ has the property $S_d(M)$ which is the key step of the whole proof, see \cref{key-prop}.
Then we show that $(x^{(i)},y,x^{(2^d-i-1)})$  also has the property $S_d(M)$ for any $i\in [2^d-1]$,
see \cref{order-111}.
Finally, based on  this property to show every point $\mathbf{x}\in A_{x_0}^{[d]}$ has the property $S_d(M_0)$
for some closed subset $M_0$ of $M$ with positive measure and $x_0\in M_0$,
we use the induction on the number of $x_0$ appearing in $\mathbf{x}$.
This is the reason why we need a closed set with  positive measure
in the definition of the property $S_d(M)$ rather than a single point.

\begin{prop}\label{order-d}
Let $M$ be a closed subset of $X$ with positive measure and $d\in\mathbb{N}$.
Then there exists a closed subset $M_0$ of $M$ with positive measure
and $x\in M_0$ such that every point in $A_x^{[d]}$ has the property $S_d(M_0)$.
\end{prop}

\begin{proof}[Proof of \cref{full-mesure=1} assuming \cref{order-d}]

For $d\in \mathbb{N}$,
let $\Theta_d$ be the metric on $X^{[d]}$ defined by
$\Theta_d(\mathbf{u},\mathbf{v})=\max_{\epsilon\in \{0,1\}^d}\rho(\mathbf{u}(\epsilon),\mathbf{v}(\epsilon))$ for $\mathbf{u},\mathbf{v}\in X^{[d]}$,
where $\rho$ is the metric on $X$.
Let $\Theta_{d,H}$ be the Hausdorff metric on $2^{X^{[d]}}$.

Consider the map
\[
F_d: X\to 2^{X^{[d]}}, \quad x\mapsto \overline{\mathcal{F}}^{[d]}(x^{[d]}).
\]
It is easy to check that $F_d$ is l.s.c..

\medskip

\noindent {\bf Claim 1:}
For any $d\in\mathbb{N}$ and any closed subset $M$ of $X$ with positive measure,
there exists $x\in M$ such that
\[
A_x^{[d]}=
\{x\}\times (\pi^{-1}(\pi(x)))^{2^d-1}\subseteq \overline{\mathcal{F}}^{[d]}(x^{[d]})=F_d(x).
\]
\begin{proof}[Proof of Claim 1]
$F_d$ is  l.s.c. and hence  Borel measurable.
By Lusin's theorem  there exists a closed subset $K$ of $X$ with $\mu(K)>1-\frac{\mu(M)}{2}$
 such that $F_d|_{K}$ is continuous.
Then we have $\mu(M\cap K)>0$.

It follows from \cref{order-d} that there exists a closed subset $M_0$ of $M\cap K$
with positive measure and $x\in M_0$ such that every point in $A_x^{[d]}$ has the property $S_d(M_0)$.
We next show that the point $x$ meets our requirement, i.e. $A_x^{[d]}\subseteq F_d(x)$.

Let $\mathbf{x}=(x_\epsilon)_{\epsilon\in\{0,1\}^d}\in A_x^{[d]}$ and fix $\gamma>0$.
As $F_d|_K$ is continuous and $x\in M_0\subseteq K$,
there exists $\delta>0$ with $\delta<\gamma$ such that whenever $x'\in K\cap B(x,\delta)$ one has $\Theta_{d,H}(F_d(x),F_d(x'))<\gamma$.
As $\mathbf{x}$ has the property $S_d(M_0)$,
there exists $w\in M_0\cap B(x,\delta)$ and $\vec{n}\in\mathbb{Z}^d$
such that
\[
w\in\bigcap_{\epsilon\in\{0,1\}^d}T^{-\vec{n}\cdot\epsilon}B(x_\epsilon,\delta).
\]

Set $\mathbf{w}=(T^{\vec{n}\cdot\epsilon}w)_{\epsilon\in\{0,1\}^d}$.
Then we have $\Theta_d(\mathbf{x},\mathbf{w})<\delta$ and $\Theta_{d,H}(F_d(w),F_d(x))<\gamma$
as $w\in M_0\cap B(x,\delta)\subseteq K\cap B(x,\delta)$.
Since $\mathbf{w}\in F_d(w)$, one can choose $\mathbf{x}'\in F_d(x)$ with $\Theta_d(\mathbf{w},\mathbf{x}')<\gamma$.
This implies that
 \[
 \Theta_d(\mathbf{x},F_d(x)) \leqslant \Theta_d(\mathbf{x},\mathbf{x}')
  \leqslant \Theta_d(\mathbf{x},\mathbf{w})+\Theta_d(\mathbf{w},\mathbf{x}')
<\delta+\gamma<2\gamma.
 \]

As $\gamma$ is arbitrary, we conclude that
$\mathbf{x}\in F_d(x)$ as was to be shown.
\end{proof}
\medskip

For $d\in\mathbb{N}$ and $\delta>0$, let
\[
\Lambda_{d,\delta}:=
\{x\in X: \exists\; \mathbf{x}\in A_x^{[d]},\;\mathrm{s.t.}\;
\Theta_d(\mathbf{x},F_d(x))\geqslant\delta
\}.
\]

\medskip

\noindent {\bf Claim 2:}
$\Lambda_{d,\delta}$ is closed.

\begin{proof}[Proof of Claim 2]
Let $\{x_i\}_{i\in\mathbb{N}}$ be a sequence in $\Lambda_{d,\delta}$ such that $\lim_{i\to\infty}x_i=x$ for some $x\in X$.
Then for every $i\in\mathbb{N}$ there exists $\mathbf{x}_{i}\in A_{x_i}^{[d]}$
such that $\Theta_d(\mathbf{x}_i,F_d(x_i))\geqslant\delta$.
Without loss of generality, we may assume that $\lim_{i\to \infty}\mathbf{x}_i=\mathbf{x}$ for some $\mathbf{x}\in X^{[d]}$.
Then $\mathbf{x}(\vec{0})=x$ and
\[
\pi^{[d]}(\mathbf{x})=\lim_{i\to\infty}\pi^{[d]}(\mathbf{x}_{i})=\lim_{i\to\infty}\pi(x_i)^{[d]}=\pi(x)^{[d]},
\]
 which implies $\mathbf{x}\in A_{x}^{[d]}$.
To show this claim, it suffices to show $\mathbf{x}\in \Lambda_{d,\delta}$.

Fix $\gamma>0$.
Recall that the map $F_d$ is l.s.c.,
there exists $\delta>0$ such that
 $F_d(x) \subseteq B(F_d(x'),\gamma/4)$ for all $x'\in B(x,\delta)$.
Choose $i\in\mathbb{N}$  such that
$\rho(x_i,x)<\delta$ and $\Theta_d(\mathbf{x}_i,\mathbf{x})<\gamma/4$.
Then $F_d(x) \subseteq B(F_d(x_i),\gamma/4)$.
Let $\mathbf{y}\in  F_d(x)$ and choose $\mathbf{y}_i\in  F_d(x_i)$ with $\Theta_d(\mathbf{y}_i,\mathbf{y})<\gamma/4$.
Then
\[
\Theta_d(\mathbf{x}_i,\mathbf{x})+\Theta_d(\mathbf{x},\mathbf{y})+\Theta_d(\mathbf{y},\mathbf{y}_i)\geqslant
\Theta_d(\mathbf{x}_i,\mathbf{y}_i)\geqslant
\Theta_d(\mathbf{x}_i,F_d(x_i))\geqslant \delta,
\]
and thus $\Theta_d(\mathbf{x},\mathbf{y})\geqslant \delta-\gamma$.
This implies  that
\[
\Theta_d(\mathbf{x},F_d(x))=\inf_{\mathbf{y}\in F_d(x)}\Theta_d(\mathbf{x},\mathbf{y})\geqslant \delta-\gamma.
\]
As $\gamma$ is arbitrary, we conclude that
$\Theta_d(\mathbf{x},F_d(x))\geqslant \delta$
and thus $\mathbf{x}\in \Lambda_{d,\delta}$ as was to be shown.
\end{proof}

It follows from Claims 1 and 2 that $\mu(\Lambda_{d,\delta})=0$ for any $d\in\mathbb{N}$ and $\delta>0$.
Now set $ \Omega= X\backslash \bigcup\limits_{d,k\in\mathbb{N}}\Lambda_{d,\frac{1}{k}}$.
Then $\mu(\Omega)=1$ and for any $x\in \Omega,d\geqslant 2$, we have
 \[
  \{x\}\times (\pi^{-1}(\pi(x)))^{2^d-1}\subseteq  \overline{\mathcal{F}^{[d]}}(x^{[d]}).
  \]

This finishes the proof.
\end{proof}

\subsection{A key proposition}\

The following proposition is a key step to the proof of \cref{order-d}.

\begin{prop}\label{key-prop}
Let $M$ be a closed subset of $X$ with positive measure.
Then there exists $x\in M$ such that
for any $\delta>0,d\in\mathbb{N}$ and any
$y\in X$ with $(x,y)\in \mathbf{RP}^{[\infty]}$, there exists a closed set $W\subseteq M\cap B(x,\delta)$
with positive measure and
$\vec{n}=(n_1,\ldots,n_d)\in\mathbb{Z}^d$ such that
\begin{itemize}
[itemsep=3pt,parsep=2pt]
  \item $T^{n_1}W\subseteq B(y,\delta)$,
  \item $T^{\sum_{i\in\epsilon}n_i}W\subseteq M\cap B(x,\delta)$ for all $\epsilon\subseteq[d],\epsilon\neq\{1\}$.
  \footnote{Here we view $\epsilon$ as a subset of $[d]$.}
  \end{itemize}
\end{prop}

To show this proposition, we need the following recurrence theorem.
We recall that a subset $R$ of $\mathbb{Z}^k$ is an {\it $IP_s$-set} if there exist $\vec{n}_1, \ldots, \vec{n}_s\in\mathbb{Z}^k$
such that \[R=
\mathrm{IP}(\vec{n}_1,\ldots,\vec{n}_s)=\{\vec{n}_{i_1}+\cdots+\vec{n}_{i_k}:
1\leqslant i_1<\cdots<i_k\leqslant s ,1\leqslant k\leqslant s\}.
\]
We shall call a subset of $\mathbb{Z}^k$ an {\it $IP^*_s$-set} if it intersects
every non-trivial $\mathrm{IP}_s$-set.

\begin{theorem}\cite[Theorem 6.15]{BM00}\label{ip-s}
  Suppose $t\in\mathbb{N}$ and $\delta$ are given.
  There exists $s=s(t,\delta)\in\mathbb{N}$ and $\xi=\xi(t,\delta)>0$ having the property that for all $r,k\in\mathbb{N}$,
  if $r$ commuting measure preserving transformations $T_1,\ldots,T_r$ of a probability space $(X,\mathcal{B},\mu)$ are given,
  as well as linear polynomials $p_{i,j}(n_1,\ldots,n_k)\in\mathbb{Z}[n_1,\ldots,n_k]$ with $p_{i,j}(0,\ldots,0)=0,1\leqslant i\leqslant r,1\leqslant j\leqslant t$, then for every $B\in\mathcal{B}$ with $\mu(B)\geqslant\delta$ the set
  \[
  R_{B,\xi}=\big\{
  (n_1,\ldots,n_k)\in\mathbb{Z}^k:
  \mu(\bigcap_{j=1}^t(\prod_{i=1}^{r}T_i^{p_{i,j}(n_1,\ldots,n_k)})^{-1}B)
  >\xi
  \big\}
  \]
  is an $IP_s^{*}$-set in $\mathbb{Z}^k$.
\end{theorem}

The following corollary follows from \cref{ip-s} directly.
\begin{cor}\label{recurrence}
Assume that $T_1,\ldots,T_d$ are commuting measure preserving transformations of a probability space $(X,\mathcal{B},\mu)$
and $B\in\mathcal{B}$ with $\mu(B)>0$,
then for any $s\in [d]$
there exists $N\in\mathbb{N}$ such that for any
$\mathbf{n}_1,\ldots,\mathbf{n}_{N}\in(\mathbb{Z}^{d})^{s}$ one can find some
$\mathbf{n}=(\vec{n}_1,\ldots,\vec{n}_s)\in \mathrm{IP}(\mathbf{n}_1,\ldots,\mathbf{n}_{N})$
such that
\[
\mu\big(
\bigcap_{\alpha\subseteq [s]}\big(\prod_{i\in \alpha}T^{-\vec{n}_i}\big)B
\big)>0,
\]
where $T^{\vec{n}}=T_1^{n_1}\cdots T_d^{n_d}$ for $\vec{n}=(n_1,\ldots,n_d)\in\mathbb{Z}$.
\end{cor}

We give some examples to explains this corollary.

Assume that $T_1,T_2$ are commuting measure preserving transformations of a probability space $(X,\mathcal{B},\mu)$
and $B\in\mathcal{B}$ with $\mu(B)>0$.
\begin{itemize}
      \item There exists $N_1\in\mathbb{N}$ such that for any
$n_1,\ldots,n_{N_1}\in \mathbb{Z}$ one can find some
$n\in \mathrm{IP}(n_1,\ldots,n_{N_1})$
such that
\[
\mu(B\cap T_1^{-n}B)>0.
\]
 \item
 There exists $N_2\in\mathbb{N}$ such that for any
$\vec{n}_1,\ldots,\vec{n}_{N_2}\in\mathbb{Z}^2$ one can find some
$\vec{n}=(a,b)\in \mathrm{IP}(\vec{n}_1,\ldots,\vec{n}_{N_2})$
such that
\[
\mu(B\cap T_1^{-a}B\cap T_2^{-b}B\cap T_1^{-a}T_2^{-b}B)>0.
\]
    \item
There exists $N_3\in\mathbb{N}$ such that for any
$\mathbf{n}_1,\ldots,\mathbf{n}_{N_3}\in(\mathbb{Z}^{2})^{2}$ one can find some
$\mathbf{n}=(\vec{u},\vec{v})\in \mathrm{IP}(\mathbf{n}_1,\ldots,\mathbf{n}_{N_3})$ with $\vec{u}=(u_1,u_2),\vec{v}=(v_1,v_2)$
such that
\begin{align*}
  &\mu(B\cap T^{-\vec{u}}B\cap T^{-\vec{v}}B\cap T^{-\vec{u}-\vec{v}}B ) \\
=&\mu(B\cap T_1^{-u_1} T_2^{-u_2}B\cap  T_1^{-v_1}T_2^{-v_2}B\cap  T_1^{-u_1-u_2}T_2^{-v_1-v_2}B)>0.
\end{align*}
\end{itemize}

\medskip

With these tools at hand we can show \cref{key-prop}.

\begin{proof}[Proof of \cref{key-prop}]
Let $M$ be a closed subset of $X$ with $\mu(M)>0$.
Let $\{\mu_x^{[d]}\}_{x\in X,d\in\mathbb{N}}$ be as in \cref{supp-order-d}.
It follows from \cref{support-d} that there exists $x\in M$ such that for any $d\in\mathbb{N}$,
\begin{enumerate}
  \item[(E)]\label{fd-ergodic} the system $(X^{[d]},\mu_x^{[d]},\mathcal{F}^{[d]})$ is ergodic;
  \item[(P)]\label{positive-measure-all-order} $\mu(M\cap B(x,\gamma))>0$ and $\mu_x^{[d]}((M\cap B(x,\gamma))^{[d]})>0$ for any $\gamma>0$.
\end{enumerate}

We claim that the point $x$ meets our requirement.

Let $y\in X$ with $(x,y)\in \mathbf{RP}^{[\infty]}$ and fix $\delta>0$.
For short, we write 
\[
V=M\cap B(x,\delta).
\]
Then by (P) we have  $\mu(V)>0$ and $\mu_x^{[d]}(V^{[d]})>0$ for any $d\in\mathbb{N}$.

\medskip

Before giving the full proof,
we will show that this proposition holds for $d=2,3$ as examples
to illustrate our basic ideas.

\medskip

\noindent {\bf The case $d=2$.}

As $\mu(V)>0$,
there exists $N\in\mathbb{N}$ such that for any $\mathrm{IP}_N$-set $R$,
one can find some $n\in R$ such that
\[
\mu(V\cap T^{-n}V)>0.
\]

By \cref{RP-infi} there is an $\mathrm{IP}_N$-set $R$ such that $T^nx\in B(y,\delta)$ for all $n\in R$.
Now we can choose $n\in R$ such that
\[
T^nx\in B(y,\delta)
\quad \text{and}\quad
\mu(V\cap T^{-n}V)>0.
\]

Choose $\eta>0$ with $\eta<\delta$ such that
\begin{equation}\label{subset2}
T^nB(x,\eta)\subseteq B(y,\delta).
\end{equation}

For short, we write $U=M\cap B(x,\eta)$.
Then $\mu(U)>0$ by (P).
Recall that $\mu$ is ergodic, there exists $m\in\mathbb{N}$ such that
\[
\mu(U\cap T^{-m}(V\cap T^{-n}V))>0.
\]
Let $W\subseteq U\cap T^{-m}V\cap T^{-n-m}V$ be a closed set with $\mu(W)>0$.
Then we have
\begin{itemize}
[itemsep=3pt,parsep=2pt]
  \item $ T^nW\subseteq T^n U\subseteq B(y,\delta)$ by \eqref{subset2},
  \item $W\subseteq U\subseteq V$ and  $T^mW,\; T^{n+m}W\subseteq V$.
\end{itemize}

This shows the case $d=2$.

\medskip

\noindent {\bf The case $d=3$.}

We first define inductively positive integers $m_3,m_2,m_1$ and $N_2,N_1$ as follows:
Set $m_3=1$.   
As $\mu(V)>0$,
it follows from \cref{recurrence} that
there exists $N_2\in\mathbb{N}$ such that for any $\vec{n}_1,\ldots,\vec{n}_{N_2}\in\mathbb{Z}^2$ one can find some
$(n,m)\in \mathrm{IP}(\vec{n}_1,\ldots,\vec{n}_{N_2})$
such that
\[
\mu(V\cap T^{-n}V\cap T^{-m}V\cap T^{-n-m}V)>0.
\]

Set $m_2=m_3N_2=N_2$.
Let us consider the system  $(X^{[N_2]},\mu_x^{[N_2]},\mathcal{F}^{[N_2]})$.
As $\mu_x^{[N_2]}(V^{[N_2]})>0$ by (P),
it follows from \cref{recurrence} that there exists $N_1\in\mathbb{N}$
such that for any $\vec{n}_1,\ldots,\vec{n}_{N_1}\in\mathbb{Z}^{m_2}=\mathbb{Z}^{N_2}$ one can find some
$\vec{n}\in \mathrm{IP}(\vec{n}_1,\ldots,\vec{n}_{N_1})$
such that $S=(T^{\vec{n}\cdot \epsilon})_{\epsilon\in\{0,1\}^{N_2}}$ and
\[
\mu_x^{[N_2]}(V^{[N_2]}\cap S^{-1}V^{[N_2]})>0.
\]

Set $m_1=m_2N_1=N_2N_1$.

\medskip

By \cref{RP-infi} there is an $\mathrm{IP}_{m_1}$-set $R$ such that $T^nx\in B(y,\delta)$ for all $n\in R$.
Choose $\eta>0$ with $\eta<\delta$ such that for all $n\in R$,
\begin{equation}\label{subset-d=3}
T^nB(x,\eta)\subseteq B(y,\delta).
\end{equation}

For short, we write $U=M\cap B(x,\eta)$.
Let $n_1^{1},\ldots,n_{m_1}^{1}\in\mathbb{N}$ such that $R=\mathrm{IP}(n_1^{1},\ldots,n_{m_1}^{1})$.
Notice that $m_1=N_2N_1$, thus
we can divide $\{n_1^{1},\ldots,n_{N_2N_1}^{1}\}$ into $N_1$ vectors of dimension $N_2$ as follows:
set $\vec{n}_i^{1}=(n_{(i-1)N_2+1}^{1},\ldots,n_{iN_2}^{1})$
for $i\in[ N_1]$.
Then 
there exists
$\vec{m}^{1}\in \mathrm{IP}(\vec{n}_1^{1},\ldots,\vec{n}_{N_1}^{1})$
such that $S_1=(T^{\vec{m}^{1}\cdot \epsilon})_{\epsilon\in\{0,1\}^{N_2}}$ and
\[
\mu_x^{[N_2]}(V^{[N_2]}\cap S_1^{-1}V^{[N_2]})>0.
\]

Recall that the system  $(X^{[N_2]},\mu_x^{[N_2]},\mathcal{F}^{[N_2]})$ is ergodic by (E)
and $\mu_x^{[N_2]}(U^{[N_2]})>0$ by (P),
there exists $\vec{m}^{2}\in\mathbb{N}^{N_2}$ such that $S_2=(T^{\vec{m}^{2}\cdot \epsilon})_{\epsilon\in\{0,1\}^{N_2}}$
and
\begin{equation}\label{positive-measure-d=3}
 \mu_x^{[N_2]}(U^{[N_2]}\cap S_2^{-1}(V^{[N_2]}\cap S_1^{-1}V^{[N_2]})) >0.
\end{equation}

Set $\vec{m}^{i}=(m_1^{i},\ldots,m_{N_2}^{i})$ for $i=1,2$ and
set $\vec{n}^{2}_j=(m_j^{1},m_j^{2})$ for $j\in[ N_2]$.
Then by the construction of $N_2$
 there exists $(a,b)\in \mathrm{IP}(\vec{n}^{2}_1,\ldots,\vec{n}^{2}_{N_2})$ such that
\[
\mu(V\cap T^{-a} V\cap T^{-b} V\cap T^{-a-b}V)>0.
\]

Let $\omega\in \{0,1\}_*^{N_2} $
such that $(a,b)=\sum_{j=1}^{N_2}\vec{n}_j^2\cdot\omega_j=(\vec{m}^{1}\cdot\omega,\vec{m}^{2}\cdot\omega)$.
Notice the $\vec{0}$-component and $\omega$-component of
 $U^{[N_2]}\cap S_2^{-1}(V^{[N_2]}\cap S_1^{-1}V^{[N_2]})$ are $U$ and 
 $W_\omega=U\cap T^{-b}V\cap T^{-a-b}V$ respectively, it follows from \cref{supp-order-d} (2)
that
\[
\mu(W_\omega)=\mu_x^{[1]}(U\times W_\omega)\geqslant
\mu_x^{[N_2]}(U^{[N_2]}\cap S_2^{-1}(V^{[N_2]}\cap S_1^{-1}V^{[N_2]}))>0.
\]

As $\mu$ is ergodic,
there exists $c\in\mathbb{N}$ such that
\[
\mu(W_\omega\cap T^{-c}(V\cap T^{-a} V\cap T^{-b} V\cap T^{-a-b}V))>0.
\]

Now let
\[
W=
U\cap T^{-b}V\cap T^{-a-b}V\cap T^{-c}V\cap T^{-c-a} V\cap T^{-c-b} V\cap T^{-c-a-b}V.
\]
Then $W\subseteq U\subseteq V$ and $T^nW\subseteq V$ for any $n\in \mathrm{IP}(a,b,c)\setminus \{a\}$.

Recall that $\vec{m}^{1}\in \mathrm{IP}(\vec{n}_1^{1},\ldots,\vec{n}_{N_1}^{1})$,
there exist $1\leqslant j_1<\cdots<j_l\leqslant N_1$ such that
$\vec{m}^{1}=\vec{n}_{j_1}^{1}+\cdots+\vec{n}_{j_l}^{1}$ which implies
\[
a=\vec{m}^1\cdot \omega=\sum_{s=1}^{l}\sum_{t\in \omega}n_{(j_s-1)N_2+t}^1
\in \mathrm{IP}(n_1^{1},\ldots,n_{m_1}^{1})=R.
\]
By \eqref{subset-d=3} we have  $T^aW\subseteq T^aU\subseteq T^aB(x,\eta)\subseteq B(y,\delta)$.

This shows the case $d=3$.

\medskip

\noindent {\bf The general case.}

We first define inductively positive integers $m_d,\ldots,m_1$ and $N_{d-1},\ldots,N_1$ as follows:
Set $m_d=1$. Assume that we have already defined $m_d,\ldots,m_{i+1}$ and $N_{d-1},\ldots,N_{i+1}$.
Let us consider the system $(X^{[m_{i+1}]},\mu_x^{[m_{i+1}]},\mathcal{F}^{[m_{i+1}]})$.
As $\mu_x^{[m_{i+1}]}(V^{[m_{i+1}]})>0$ by (P),
it follows from \cref{recurrence} that there exists $N_i\in\mathbb{N}$ such that for any
$\mathbf{n}_1,\ldots,\mathbf{n}_{N_i}\in(\mathbb{Z}^{m_{i+1}})^{i}$ one can find some
$\mathbf{n}=(\vec{n}_1,\ldots,\vec{n}_{i})\in \mathrm{IP}(\mathbf{n}_1,\ldots,\mathbf{n}_{N_i})$
such that $S_j=(T^{\vec{n}_j\cdot \epsilon})_{\epsilon\in\{0,1\}^{m_{i+1}}},j\in [i]$ and
\[
\mu_x^{[m_{i+1}]}
\Big(\bigcap_{\alpha\subseteq [i]}\big(\prod_{j\in\alpha}S_{j}^{-1}\big) V^{[m_{i+1}]}
\Big)>0.
\]
Set $m_{i}=m_{i+1}N_i$.
This finishes the inductive definition.

\medskip

By \cref{RP-infi} there is an $\mathrm{IP}_{m_1}$-set $R$ such that $T^nx\in B(y,\delta)$ for all $n\in R$.
Choose $\eta>0$ with $\eta<\delta$ such that for all $n\in R$,
\begin{equation}\label{subset-d-general}
  T^nB(x,\eta)\subseteq B(y,\delta).
\end{equation}

For $i\in\{2,\ldots,d\}$,
inductively we will construct
$\vec{n}_1^{i},\ldots,\vec{n}_i^{i}\in\mathbb{N}^{m_i}$ and sets $W^i_\epsilon\subseteq M\cap B(x,\eta),\epsilon\in\{0,1\}^{m_i}$
with $W_{\vec{0}}^i=M\cap B(x,\eta)$
such that
\begin{enumerate}[itemsep=4pt,parsep=2pt,label=(\arabic*)]
\item[$(A_i)$]  $\vec{n}_1^i\cdot\epsilon\in R$, i.e. $  T^{\vec{n}_1^i\cdot\epsilon}B(x,\eta)\subseteq B(y,\delta)$
   for any $\epsilon\in\{0,1\}^{m_i}_*$,
\item[$(B_i)$] $\mu_x^{[m_i]}\big(\prod\limits_{\epsilon\in\{0,1\}^{m_i}}W^i_\epsilon\big)>0$,
\item[$(C_i)$]  $\big( \prod\limits_{j\in\alpha}T^{\vec{n}_{j}^{i}\cdot \epsilon}\big)W_{\epsilon}^{i}\subseteq
V$  for any $\epsilon\in\{0,1\}^{m_i}$ and any $\alpha\subseteq [i],  \alpha\neq \{1\}$.
\end{enumerate}

Assume this has been achieved, then
there exist
$n_1^d,\ldots,n_d^d\in\mathbb{N}^{m_d}=\mathbb{N}$ and sets
$W_0^d,W_1^d\subseteq M\cap B(x,\eta)$ such that
\begin{itemize}[itemsep=3pt,parsep=2pt]
\item  $T^{n_1^d}B(x,\eta)\subseteq B(y,\delta)$,
\item $\mu(W_1^d)=\mu_x^{[1]}(W_0^d\times W_1^d)>0$,
\item $T^{\sum_{j\in\alpha}n_j^d}W_1^d=
\big(\prod_{j\in\alpha}T^{n_j^d}\big)W_1^d\subseteq V$ for any $\alpha\subseteq [i],  \alpha\neq \{1\}$.
\end{itemize}
Let $W\subseteq W_1^d$ be a closed set with positive measure and set $n_i=n_i^d$ for $i\in [d]$.
Then we have
\begin{itemize}[itemsep=3pt,parsep=2pt]
  \item $T^{n_1}W\subseteq T^{n_1^d}W_1^d\subseteq T^{n_1^d}B(x,\eta)\subseteq B(y,\delta)$,
        \item $T^{\sum_{j\in \alpha}n_j}W\subseteq T^{\sum_{j\in \alpha}n_j^d}W_1^d\subseteq V=M\cap B(x,\delta)$
        for any $\alpha\subseteq [i],  \alpha\neq \{1\}$,
\end{itemize}
which shows the proposition.

\medskip

{\bf We now return to the inductive construction.}

\medskip

\noindent {\bf Step 1:}
For short, we write $U=M\cap B(x,\eta)$.
Let $n_1^{1},\ldots,n_{m_1}^{1}\in\mathbb{N}$ such that $R=\mathrm{IP}(n_1^{1},\ldots,n_{m_1}^{1})$.
Since $m_1=m_{2}N_1$,
we divide $\{n_1^{1},\ldots,n_{m_1}^{1}\}$ into $N_1$ vectors of dimension $m_2$ as follows:
set $\vec{r}_{k}=(n_{(k-1)m_{2}+1}^1,\ldots,n_{km_{2}}^1)$
for $ k\in[N_1]$.

Let us consider the system $(X^{[m_{2}]},\mu_x^{[m_{2}]},\mathcal{F}^{[m_{2}]})$.
Recall that $\mu_x^{[m_{2}]}(V^{[m_2]})>0$ by (P),
thus by the construction of $N_1$ there exists a non-empty subset $\beta$ of $ [N_1]$ such that
\[
\vec{n}_1^{2}=\sum_{k\in\beta}\vec{r}_{k},\quad
S_1^2=(T^{-\vec{n}_1^{2}\cdot \epsilon})_{\epsilon\in\{0,1\}^{m_{2}}},
\]
and
\[
\mu_x^{[m_2]}(V^{[m_2]}\cap S_1^2 V^{[m_2]})>0.
\]

As $\mu_x^{[m_{2}]}$ is $\mathcal{F}^{[m_{2}]}$-ergodic by (E) and $\mu_x^{[m_2]}(U^{[m_2]})>0$ by (P), there exists $\vec{n}_{2}^{2}\in\mathbb{N}^{m_{2}}$
such that $S_2^{2}=(T^{-\vec{n}_{2}^{2}\cdot \epsilon})_{\epsilon\in\{0,1\}^{m_{2}}}$ and
\[
\mu_x^{[m_{2}]}
(U^{[m_2]}\cap S_2^2(V^{[m_2]}\cap S_1^2 V^{[m_2]}))
>0.
\]

Set $W_\epsilon^2=U\cap T^{-\vec{n}_2^2\cdot\epsilon}V\cap T^{-(\vec{n}_2^2+\vec{n}_1^2)\cdot\epsilon}V\subseteq U$ for $\epsilon\in\{0,1\}^{m_2}$.
Then we have
\begin{itemize}[itemsep=3pt,parsep=2pt]
  \item for any $\epsilon\in\{0,1\}^{m_2}_*$,
  \[
\vec{n}_1^{2}\cdot \epsilon=\sum_{k\in\beta}\vec{r}_{k}\cdot \epsilon
=\sum_{k\in\beta}\sum_{j=1}^{m_2}n_{(k-1)m_{2}+j}^1\cdot \epsilon_j\in R,
\]
  \item $\mu_x^{[m_2]}\big(\prod\limits_{\epsilon\in\{0,1\}^{m_2}}W^2_\epsilon\big)=\mu_x^{[m_{2}]}
(U^{[m_2]}\cap S_2^2(V^{[m_2]}\cap S_1^2 V^{[m_2]}))
>0$,
 \item
 $W_{\epsilon}^2\subseteq
 V\cap T^{-\vec{n}_2^2\cdot\epsilon}V\cap T^{-(\vec{n}_2^2+\vec{n}_1^2)\cdot\epsilon}V$
 for any $\epsilon\in\{0,1\}^{m_2}$.
\end{itemize}

\medskip

\noindent {\bf Step i:}
Let $i\geqslant 2$ be an integer and assume that we have already chosen
$\vec{n}_1^{i},\ldots,\vec{n}_i^{i}\in\mathbb{N}^{m_i}$ and sets $W^i_\epsilon\subseteq U,\epsilon\in\{0,1\}^{m_i}$
such that
\begin{enumerate}[itemsep=4pt,parsep=2pt,label=(\arabic*)]
\item[$(A_i)$]  $\vec{n}_1^i\cdot\epsilon\in R$ for any $\epsilon\in\{0,1\}^{m_i}_*$,
\item[$(B_i)$] $\mu_x^{[m_i]}\big(\prod\limits_{\epsilon\in\{0,1\}^{m_i}}W^i_\epsilon\big)>0$,
\item[$(C_i)$]  $\big( \prod\limits_{j\in\alpha}T^{\vec{n}_{j}^{i}\cdot \epsilon}\big)W_{\epsilon}^{i}\subseteq
V$  for any $\epsilon\in\{0,1\}^{m_i}$ and any $\alpha\subseteq [i],  \alpha\neq \{1\}$.
\end{enumerate}



For $j\in[i]$ set $\vec{n}_j^{i}=(n_{j,1}^{i},\ldots,n_{j,m_i}^{i})$.
As $m_{i}=m_{i+1}N_i$,
we can divide $\{n_{j,1}^{i},\ldots,n_{j,m_i}^{i}\}$ into $N_{i}$ vectors of dimension $m_{i+1}$ as follows:
set $\vec{r}_{j,k}=(n_{j,(k-1)m_{i+1}+1}^{i},\ldots,n_{j,km_{i+1}}^{i})$
for $ k\in [N_i]$.
Set $\mathbf{n}_k=(\vec{r}_{1,k},\ldots,\vec{r}_{i,k})\in(\mathbb{Z}^{m_i})^i$ for $k\in [N_i]$.

Let us consider the system $(X^{[m_{i+1}]},\mu_x^{[m_{i+1}]},\mathcal{F}^{[m_{i+1}]})$ and the set $V^{[m_{i+1}]}$.
By the construction of $N_i$ there exist $1\leqslant k_1<\cdots<k_l\leqslant N_i$ such that for $j\in [i]$
\[
\vec{n}_j^{i+1}=\vec{r}_{j,k_1}+\cdots+\vec{r}_{j,k_l},\quad
S_{j}^{i+1}=(T^{-\vec{n}_{j}^{i+1}\cdot \epsilon})_{\epsilon\in\{0,1\}^{m_{i+1}}},
\]
and
 \[
 \mu_x^{[m_{i+1}]}\Big(
 \underbrace{\bigcap_{ \alpha\subseteq [i]   }
 \big(\prod_{j\in\alpha}S_{j}^{i+1}\big)V^{[m_{i+1}]}}_{=:\mathbf{V}_1}
 \Big)>0.
  \]

For $s\in [m_{i+1}]$, set
\[
\omega_s=\{(k_1-1)m_{i+1}+s,\ldots,(k_l-1)m_{i+1}+s\}.
\]
Then $\omega_1,\ldots,\omega_{m_{i+1}}$ are non-empty pairwise disjoint subsets of $[m_i]$.
Let  $p:\{0,1\}^{m_{i+1}}\to \{0,1\}^{m_i}$ be the map such that $p(\emptyset)=\emptyset$ and
$p(\eta)=\bigcup_{s\in\eta} \omega_s
$ for $\eta\in \{0,1\}^{m_{i+1}}_*$. 
It follows from \cref{supp-order-d} (2) that
\[
 \mu_x^{[m_{i+1}]}\Big(
\underbrace{\prod_{\eta\in\{0,1\}^{m_{i+1}}}W^i_{p(\eta)}}_{=:\mathbf{V}_2}
\Big)
\geqslant
 \mu_x^{[m_{i}]}\Big(
\prod_{\epsilon\in\{0,1\}^{m_{i}}}W^i_{\epsilon}
\Big)
>0.
\]

Moreover  for $j\in [i]$ and $\eta\in \{0,1\}^{m_{i+1}}_*$, we have
\begin{equation}\label{low-equation}
\begin{split}
\vec{n}_j^{i+1}\cdot\eta  =\sum_{t=1}^{l}\vec{r}_{j,k_t}\cdot\eta&=\sum_{t=1}^{l}\sum_{s\in\eta}n_{j,(k_t-1)m_{i+1}+s} \\
     & =\sum_{s\in\eta}\vec{n}_j^i\cdot\omega_s=
\vec{n}_j^i\cdot p(\eta).
\end{split}
\end{equation}
Recall that $\mu_x^{[m_{i+1}]}$ is $\mathcal{F}^{[m_{i+1}]}$-ergodic by (E),
there exists $\vec{n}_{i+1}^{i+1}\in\mathbb{N}^{m_{i+1}}$ such that
$S_{i+1}^{i+1}=(T^{-\vec{n}_{i+1}^{i+1}\cdot \eta})_{\eta\in\{0,1\}^{m_{i+1}}}$
and
\begin{equation}\label{positive-A}
  \mu_x^{[m_{i+1}]}
(\mathbf{V}_1\cap S_{i+1}^{i+1}\mathbf{V}_2)
>0.
\end{equation}

For $\eta\in\{0,1\}^{m_{i+1}}$ set
\[
W_{\eta}^{i+1}=W^i_{p(\eta)}
\cap T^{-\vec{n}_{i+1}^{i+1}\cdot \eta}\Big(
\bigcap_{ \alpha\subseteq [i]   }
 \big(\prod_{j\in\alpha}T^{-\vec{n}_{j}^{i+1}\cdot \eta}\big)V\Big).
\]
Then we have
\begin{align}
\label{A1}\prod_{\eta\in\{0,1\}^{m_{i+1}}}W_{\eta}^{i+1}& =\mathbf{V}_1\cap S_{i+1}^{i+1}\mathbf{V}_2, \\
\label{A2} T^{\vec{n}_{i+1}^{i+1}\cdot \eta}
 \big(\prod_{j\in\alpha}T^{\vec{n}_{j}^{i+1}\cdot \eta}\big) W_{\eta}^{i+1} &\subseteq V,\quad \forall\; \alpha\subseteq [i],\\
 \label{A3}  W_{\eta}^{i+1}&\subseteq W^i_{p(\eta)}.
\end{align}

We next show that $\vec{n}_1^{i+1},\ldots,\vec{n}_{i+1}^{i+1}$ and sets $W^{i+1}_\eta,\eta\in\{0,1\}^{m_{i+1}}$ satisfy $(A)_{i+1},(B)_{i+1},(C)_{i+1}$.

By \eqref{low-equation} and $(A)_i$ for any $\eta\in\{0,1\}_*^{m_{i+1}}$ we have
\[
\vec{n}_1^{i+1}\cdot\eta=\vec{n}_1^i\cdot p(\eta)\in R.
\]
By \eqref{positive-A} and \eqref{A1} we have
\[
\mu_x^{[m_{i+1}]}\big(\prod\limits_{\eta\in\{0,1\}^{m_{i+1}}}W^{i+1}_\eta\big)=\mu_x^{[m_{i+1}]}(\mathbf{V}_1\cap S_{i+1}^{i+1}\mathbf{V}_2)>0.
\]

It remains to show $(C)_{i+1}$.
 Let $\alpha\subseteq [i+1]$ with $\alpha\neq \{1\}$.

When $i+1\in \alpha$. Set $\alpha'=\alpha\backslash\{i+1\}$. Then by \eqref{A2} we have
 \[
 \big( \prod\limits_{j\in\alpha}T^{\vec{n}_{j}^{i+1}\cdot \eta}\big)W_{\eta}^{i+1}=
 T^{\vec{n}_{i+1}^{i+1}\cdot \eta}
  \big( \prod\limits_{j\in\alpha'}T^{\vec{n}_{j}^{i+1}\cdot \eta}\big)W_{\eta}^{i+1}
 \subseteq
V.
\]

When $i+1\not\in \alpha$, we can view $\alpha$ as a subset of $[i]$.
Hence by $(C)_i$ \eqref{low-equation} and \eqref{A3} we have
 \[
 \big(\prod_{j\in\alpha}T^{\vec{n}_{j}^{i+1}\cdot \eta}\big) W_\eta^{i+1}=
  \big( \prod_{j\in\alpha}T^{\vec{n}_{j}^{i}\cdot p(\eta)} \big) W_\eta^{i+1}
  \subseteq  \big( \prod_{j\in\alpha}T^{\vec{n}_{j}^{i}\cdot p(\eta)} \big) W_{p(\eta)}^{i}
  \subseteq V.
 \]

We finish the construction by induction and hence ends the proof.
\end{proof}

\subsection{Proof of \cref{order-d}}\

For $x\in X,d\in \mathbb{N}$ and $i\in\{0,1,\ldots,2^d-1\}$, let
\[
D_x^{[d]}(i)  =
\{ \mathbf{x}\in A_x^{[d]}:
 | \{\epsilon\in\{0,1\}^d
 :\mathbf{x}(\epsilon)\neq x\}| =i
 \}.
\]
Then $D_x^{[d]}(0)=\{x^{[d]}\}$ and $A_x^{[d]}=\bigcup_{i=0}^{2^d-1}D_x^{[d]}(i)$.

\begin{prop}\label{order-111}
Let $M$ be a closed subset of $X$ with positive measure and $d\in \mathbb{N}$.
Then there exists $x\in M$ such that every point in $D_x^{[d]}(1)$ has the property $S_d(M)$.
\end{prop}

\begin{proof}
Let $x\in M$ be as in \cref{key-prop} for $M$.
We claim that $x$ meets our requirement.
That is, every point in $D_x^{[d]}(1)$ has the property $S_d(M)$ .

Let $\mathbf{x}\in D_x^{[d]}(1)$ and
let $\omega\in \{0,1\}^d$ such that $\mathbf{x}(\omega)=y\neq x$.
Fix $\delta>0$.
It follows from \cref{key-prop} that
there exists a closed set $W_x\subseteq M\cap B(x,\delta)$
with positive measure and
$\vec{m}=(m_1,\ldots,m_d)\in\mathbb{Z}^d$ such that
\begin{itemize}
[itemsep=3pt,parsep=2pt]
  \item $T^{m_1}W_x\subseteq B(y,\delta)$,
    \item $T^{\sum_{i\in\epsilon}m_i}W\subseteq M\cap B(x,\delta)$ for all $\epsilon\subseteq[d],\epsilon\neq\{1\}$.
\end{itemize}

When $\omega=\{1\}$, by taking $W=W_x$ and $n_i=m_i,i\in[d]$
the proposition follows.

When $\omega\neq\{1\}$,
set $j=\min \omega$ and
let $\vec{n}=(n_1,\ldots,n_d)\in\mathbb{Z}^d$ such that
\[
n_i=
\begin{cases}
m_j, &i=1, \\
m_1,& i=j,\\
-m_i,& i\in \omega\backslash\{j\},\\
m_i,& \mathrm{otherwise},
\end{cases}
\]
and $W=T^{\sum_{i\in  \omega\backslash\{j\}}m_i}W_x$.
Then $W$ is a closed subset of $X$ with positive measure.

Notice that
$\vec{n}\cdot\omega+\sum_{i\in  \omega\backslash\{j\}}m_i=m_1$
and
\[
\{\vec{n}\cdot\epsilon+\sum_{i\in  \omega\backslash\{j\}}m_i:\epsilon\in\{0,1\}^d\}
=\{\vec{m}\cdot\epsilon:\epsilon\in\{0,1\}^d\},
\]
we get
\begin{itemize}
[itemsep=3pt,parsep=2pt]
  \item $T^{\vec{n}\cdot \omega}W=T^{\vec{n}\cdot\omega+\sum_{i\in  J\backslash\{j\}}m_i}W_x
  =T^{m_1}W_x\subseteq B(y,\delta)$,
  \item $T^{\vec{n}\cdot \epsilon}W=T^{\vec{n}\cdot\epsilon+\sum_{i\in  J\backslash\{j\}}m_i}W_x\subseteq M\cap B(x,\delta)$ for all $\epsilon\in\{0,1\}^d\backslash\{\omega\}$.
\end{itemize}
In particular, $W=T^{\vec{n}\cdot \vec{0}}W\subseteq M\cap B(x,\delta)$.

This finishes the proof.
\end{proof}

Now we can give a proof of \cref{order-d}.

\begin{proof}[Proof of \cref{order-d}]
Let $M$ be a closed subset of $X$ with positive measure and $d\in\mathbb{N}$.
Consider the map
\[
\Gamma: X\to 2^{X}, \ x\mapsto \pi^{-1}(\pi(x)).
\]
Then $\Gamma$ is u.s.c. and hence Borel measurable. By Lusin's theorem, there exists a closed subset $M_0$ of $M$ with positive measure
such that $\Gamma|_{M_0}$ is continuous.

Let $x\in M_0$ be as in \cref{order-111} for $M_0$.
We claim that $x$ meets our requirement.
To prove this claim,
we will show by induction on $i$ that every point in $D_x^{[d]}(i)$ has the property $S_d(M_0)$.

Fix $\delta>0$. As the map $\Gamma$ is continuous on $M_0$,
there exists $\gamma>0$ with $\gamma<\delta$ such that
whenever $w\in M_0\cap B(x,\gamma)$ one has $\rho_H(\pi^{-1}(\pi(w)),\pi^{-1}(\pi(x)))<\delta$,
where $\rho_H$ is the Hausdorff metric on the hyperspace of $X$.
That is, for any $w\in M_0\cap B(x,\gamma)$ and any $x'\in X$ with $(x,x')\in \mathbf{RP}^{[\infty]}$
there exists $w'\in X$ with $(w,w')\in \mathbf{RP}^{[\infty]}$
such that $\rho(w',x')<\delta$.

\medskip

Before the induction we show a special case when $d=2$ and $i=2$ as an example, i.e., we will show that
$(x,y,z,x)$ has the property $S_2(M_0)$  if $\pi(y)=\pi(z)=\pi(x)$.

\medskip

\noindent {\bf A special case.}

By \cref{order-111} $(x,x,x,z)$ has the property $S_2(M_0)$.
Then
there exists a closed subset $W_z\subseteq M_0\cap B(x,\gamma)$ with positive measure and $(a_1,a_2)\in\mathbb{Z}^2$ such that
\begin{itemize}[itemsep=3pt,parsep=2pt]
\item $T^{a_1+a_2}W_z\subseteq B(z,\gamma)$,
\item $T^{a_1}W_z,T^{a_2}W_z\subseteq M_0\cap B(x,\gamma)$.
\end{itemize}

As $\mu(W_z)>0$, it follows from \cref{order-111} that
there exists $w\in W_z$ such that every point in $D_w^{[d]}(1)$ has the property $S_2(W_z)$.
As $w\in W_z\subseteq M_0\cap B(x,\gamma)$ and $(x,y)\in \mathbf{RP}^{[\infty]}$, we can choose $w'\in B(y,\delta)$ with $(w,w')\in \mathbf{RP}^{[\infty]}$.
Then $(w,w',w,w)$ has the property $S_2(W_z)$.
Choose $\eta>0$ with $B(w',\eta)\subseteq B(y,\delta)$.
Hence there exists a closed subset $W_w\subseteq W_z\cap B(w,\eta)$ with positive measure and $(b_1,b_2)\in\mathbb{Z}^2$ such that
\begin{itemize}[itemsep=3pt,parsep=2pt]
\item $T^{b_1}W_w\subseteq B(w',\eta)$,
\item $T^{b_2}W_w,T^{b_1+b_2}W_w\subseteq W_z\cap B(w,\eta)$.
\end{itemize}
Set $n_1=b_1-a_1,n_2=b_2+a_2$ and $W=T^{a_1}W_w$.
Then $\mu(W)=\mu(W_w)>0$ and
\begin{itemize}[itemsep=3pt,parsep=2pt]
\item $W=T^{a_1}W_w\subseteq T^{a_1}W_z\subseteq M_0\cap B(x,\delta)$,
\item $T^{n_1}W=T^{b_1-a_1}(T^{a_1}W_w)=T^{b_1}W_w\subseteq B(w',\eta)\subseteq B(y,\delta)$,
\item $T^{n_2}W=T^{b_2+a_2}(T^{a_1}W_w)=T^{a_1+a_2}(T^{b_2}W_w)
\subseteq T^{a_1+a_2}W_z\subseteq B(z,\delta)$,
\item $T^{n_1+n_2}W=
T^{b_1-a_1+b_2+a_2}(T^{a_1}W_w)=T^{a_2}(T^{b_1+b_2}W_w)\subseteq T^{a_2}W_z
\subseteq M_0\cap B(x,\delta)$.
\end{itemize}

This shows that $(x,y,z,x)$ has the property $S_2(M_0)$.

\medskip

\noindent {\bf The general case.}

When $i=1$, it follows from \cref{order-111}.

Let $i\geqslant 2$ be an integer and
assume that we have already shown that every point in $D_x^{[d]}(j)$ has the property $S_d(M_0)$ for all $j\in[i-1]$.

We next show that every point in $D_x^{[d]}(i)$ also has the property $S_d(M_0)$.
Let $(x_\epsilon)_{\epsilon\in \{0,1\}^d}\in D_x^{[d]}(i)$ and let
\[
\Lambda=\{\epsilon\in \{0,1\}^d:x_\epsilon\neq x\}.
\]
Then $|\Lambda|=i$. Choose $\omega\in \Lambda$ and
let $\Phi$ be a permutation  of $\{0,1\}^d$ such that for any $\epsilon\in \{0,1\}^d$,
\[
\Phi(\epsilon)_i=
\begin{cases}
1-\epsilon_i,&\mathrm{ if} \; i\in \omega,\\
\epsilon_i,& \mathrm{otherwise}.
\end{cases}
\]
Then we have $\Phi(\vec{0})=\omega,\Phi(\omega)=\vec{0}$ and $\Phi(\Phi(\epsilon))=\epsilon$ for any $\epsilon\in \{0,1\}^d$.

Let $(y_\epsilon)_{\epsilon\in\{0,1\}^d}\in X^{[d]}$ such that $y_{\vec{0}}=x$ and $y_\epsilon=x_{\Phi(\epsilon)}$ otherwise.
Then we have
\[
\Lambda'=
\{\epsilon\in\{0,1\}^d:y_\epsilon\neq x\}
=\Phi(\Lambda)\backslash\{\vec{0}\},
\]
and thus $(y_\epsilon)_{\epsilon\in\{0,1\}^d}\in D_{x}^{[d]}(i-1)$.
By our inductive assumption, the point $(y_\epsilon)_{\epsilon\in\{0,1\}^d}$ has the property $S_d(M_0)$.
Then
there exists a closed subset $W_1\subseteq M_0\cap B(x,\gamma)$ with positive measure and $\vec{a}\in\mathbb{Z}^d$ such that
\begin{equation}\label{general-a}
  W_1\subseteq \bigcap_{\epsilon\in\{0,1\}^d}T^{-\vec{a}\cdot \epsilon} B_{M_0}(y_\epsilon,\gamma),
\end{equation}
where
\[
B_{M_0}(y_\epsilon,\gamma)=
\begin{cases}
M_0\cap B(x,\gamma), &\mathrm{if} \; y_\epsilon=x, \\
 B(y_\epsilon,\gamma),& \mathrm{otherwise}.
\end{cases}
\]

On the other hand, using \cref{key-prop} for the set $W_1$,
there exists $z\in W_1$ such that every point in $D_z^{[d]}(1)$ has the property $S_d(W_1)$.
As $z\in W_1\subseteq M_0\cap B(x,\gamma)$, we can choose $z'\in B(x_\omega,\delta)$ with $(z,z')\in \mathbf{RP}^{[\infty]}$.
Choose $\eta>0$ with $B(z',\eta)\subseteq B(x_\omega,\delta)$.
Then there exists a closed subset $W_2\subseteq W_1\cap B(z,\eta)$ with positive measure and $\vec{b}\in\mathbb{Z}^d$ such that
\begin{equation}\label{general-b}
  W_2\subseteq \bigcap_{\epsilon\in\{0,1\}^d}T^{-\vec{b}\cdot \epsilon} B_{W_1}(z_\epsilon,\eta),
\end{equation}
where
\[
B_{W_1}(z_\epsilon,\eta)=
\begin{cases}
B(z',\eta), &\mathrm{if} \; \epsilon=\omega, \\
W_1\cap  B(x,\eta),& \mathrm{otherwise}.
\end{cases}
\]

Let $\vec{c}=(c_1,\ldots,c_d)\in\mathbb{Z}^d$ such that
\[
c_i=
\begin{cases}
-a_i,&\mathrm{ if} \; i\in \omega,\\
a_i,& \mathrm{otherwise}.
\end{cases}
\]

Set $\vec{n}=\vec{b}+\vec{c}$ and $W=T^{\vec{a}\cdot \omega}W_2$.
Then $W$ is a closed subset of $X$ with positive measure.
We next show
\[
  W\subseteq \bigcap_{\epsilon\in\{0,1\}^d}T^{-\vec{n}\cdot \epsilon} B_{M_0}(x_\epsilon,\delta),
\]
where
\[
B_{M_0}(x_\epsilon,\delta)=
\begin{cases}
M_0\cap B(x_\epsilon,\delta), &\mathrm{if} \; x_\epsilon=x, \\
B(x_\epsilon,\delta) ,& \mathrm{otherwise}.
\end{cases}
\]

Notice that for  $\epsilon\in \{0,1\}^d$,
\begin{equation*}
\begin{aligned}
& (\vec{b}+\vec{c})\cdot\epsilon+\vec{a}\cdot \omega\\
=\;& \vec{b}\cdot\epsilon+\sum_{i\in [d]\backslash \omega}a_i\epsilon_i-\sum_{i\in \omega}a_i\epsilon_i+\sum_{i\in \omega}a_i\\
=\; &\vec{b}\cdot\epsilon+\sum_{i\in [d]\backslash \omega}a_i\epsilon_i+\sum_{i\in\omega}a_i(1-\epsilon_i)\\
=\; &\vec{b}\cdot\epsilon+\vec{a}\cdot \Phi(\epsilon),
\end{aligned}
\end{equation*}
thus we have
\begin{equation}\label{eq:subset}
T^{\vec{n}\cdot\epsilon}W=T^{(\vec{b}+\vec{c})\cdot\epsilon}(T^{\vec{a}\cdot \epsilon}W_2)=
T^{\vec{a}\cdot \Phi(\epsilon)}(T^{\vec{b}\cdot\epsilon}W_2).
\end{equation}

When $\epsilon=\omega$, as $ \Phi(\omega)=\vec{0}$, by \eqref{general-b} and \eqref{eq:subset} we have
\[
T^{\vec{n}\cdot\omega}W=T^{\vec{b}\cdot\omega}W_2\subseteq B_{W_1}(z_{\omega},\eta)=
B(z',\eta)\subseteq B(x_\omega,\delta).
\]

When $\epsilon\in \{0,1\}^d\backslash\{\omega\}$,
 by \eqref{eq:subset}
we have
\begin{align*}
T^{\vec{n}\cdot\epsilon}W &=T^{\vec{a}\cdot \Phi(\epsilon)}(T^{\vec{b}\cdot\epsilon}W_2)&\\
&\subseteq  T^{\vec{a}\cdot \Phi(\epsilon)} (W_1\cap B(x,\eta))& \text{by} \;(\ref{general-b})\\
 & \subseteq B_{M_0}(y_{\Phi(\epsilon)},\gamma)&\text{by} \;(\ref{general-a}) \\
 &=B_{M_0}(x_{\Phi(\Phi(\epsilon))},\gamma)\subseteq
  B_{M_0}(x_{\epsilon},\delta) .&  ( \mathrm{as} \;\Phi(\Phi(\epsilon))=\epsilon)
  \end{align*}
In particular,
$W=T^{\vec{n}\cdot \vec{0}}W\subseteq  B_{M_0}(x_{\vec{0}},\delta)= B_{M_0}(x,\delta)$.

This completes the induction
procedure and hence the proof of the proposition.
\end{proof}

\section{Proof of \cref{regular-d}}

To end this paper, using \cref{full-mesure=1} we show \cref{regular-d}.

\begin{proof}[Proof of \cref{regular-d}]
Let $(X,T)$ be a minimal system without any non-trivial $\mathrm{Ind}_{fip}$-$(k+1)$-tuples.
Let $\pi:X\to X/\mathbf{RP}^{[\infty]}=Y$ be the factor map
and suppose that $m\in M(Y,T)$ is the unique ergodic measure.

For $j\in\mathbb{N}\cup\{\infty\}$, set
\[
 Y_j  =\{y\in Y:  |\pi^{-1}(y)|=j\}.
\]
By definition, $\{Y_j\}_{j\in\mathbb{N}\cup\{\infty\}}$ is a disjoint family of $Y$.
It was proved in \cite[Proposition 3.7]{HLSY21} that
$Y_j$ is measurable and $T$-invariant  for every $j\in\mathbb{N}\cup\{\infty\}$.
Since $m$ is ergodic, there exists $j\in\mathbb{N}\cup\{\infty\}$ such that $m(Y_{j})=1$ and $m(Y_i)=0$ for all $i\neq j$.

\medskip

\noindent {\bf Claim 1:}  $j\leqslant k$.

\begin{proof}[Proof of Claim 1]
Suppose for a contradiction that $j\geqslant k+1$.
Let $\mu\in M(X,T)$ be an ergodic measure.
It follows from \cref{full-mesure=1} that there exists a subset $\Omega$ of $ X$ with $\mu(\Omega)=1$ such that
$\{x\}\times (\pi^{-1}(\pi(x)))^{2^d-1}\subseteq  \overline{\mathcal{F}^{[d]}}(x^{[d]})$  for any $x\in \Omega$ and any $d\in\mathbb{N}$.
Notice that $(Y,T)$ is uniquely ergodic, we have $\pi_*(\mu)=m$ and thus $\mu(\Omega\cap \pi^{-1}(Y_{j}))=1$.
Now let $x\in \Omega\cap \pi^{-1}(Y_{j})$ and let $x_1,\ldots,x_{k+1}\in \pi^{-1}(x)$ be distinct points.
It follows that $\{x\}\times \{x_1,\ldots,x_{k+1}\}^{2^{d+1}-1}\subseteq \overline{\mathcal{F}^{[d+1]}}(x^{[d+1]})$ for any $d\geqslant 1$
and thus
$\{x_1,\ldots,x_{k+1}\}^{[d]}\subseteq \mathbf{Q}^{[d]}(X)$ for any $d\geqslant 1$ which implies
$(x_1,\ldots,x_{k+1})$ is an $\mathrm{Ind}_{fip}$-$(k+1)$-tuple by \cref{tuple}, a contradiction.
From this, we conclude $j\leqslant k$, which shows Claim 1.
\end{proof}

\noindent {\bf Claim 2:}  $X$ has at most $j$ ergodic measures.

\begin{proof}[Proof of Claim 2]
Suppose for a contradiction that $\mu_1,\ldots,\mu_{j+1}$ are distinct ergodic measures of $(X,T)$.
Let $W_i$ be the set of all $\mu_i$-generic points for $ i\in [j+1]$.
By \cref{generic-points} we have $\mu_i(W_i)=1$ for every $ i\in [j+1]$.
By Lusin's theorem, $\{\pi(W_i)\}_{i=1}^{j+1}$ are universally measurable as they are analytic sets.
\footnote{Let $X$ be a standard Borel space and let $\mathcal{B}(X)$ be the set of Borel sets of $X$. A subset
$A$ of $X$ is analytic if there exist a standard Borel space $X', A'\in \mathcal{B}(X')$ and a measurable
map $f : X'\to X$ such that $A = f(A')$. Let $\mu$ be a Borel measure on $X$. A subset $A$ is
called $\mu$-measurable set if $A = B \Delta N$ for some $B \in \mathcal{B}(X)$ and $\mu$-null set $N$. A subset
$A\subseteq  X$ is called universally measurable if it is $\mu$-measurable for every Borel measure $\mu$.
Lusin Theorem says that every analytic set is universally measurable (see \cite[Theorem 21.10]{KAS95} for example).
}
Notice that $\pi_*(\mu_1)=\cdots=\pi_*(\mu_{j+1})=m$, hence we have
\[
m(\pi(W_1))=\cdots=m(\pi(W_{j+1}))=1.
\]
Set $Y'=\bigcap_{i=1}^{j+1}\pi(W_i)$.
It follows that $m(Y')=1$ and $|\pi^{-1}(y)|\geqslant j+1$ for any $y\in Y'$.
By Claim 1 we have $m(Y_j)=1$ and thus $Y'\cap Y_j \neq\emptyset$, a contradiction.

This shows Claim 2.
\end{proof}

It remains to show that
$\pi$ is an almost $j'$ to one extension for some $j'\in[j]$.
As a matter of fact,
it follows from Claim 1 that $Y_j$ is non-empty.
Then by \cref{almost-finite-to-one} there exists $j'\in[j]$ such that $Y_{j'}$ contains a dense $G_\delta$ set, i.e., $\pi$ is an almost $j'$ to one extension.

This finishes the proof.
\end{proof}


\begin{thebibliography}{SSS}

\bibitem{BM00}
V. Bergelson and R. McCutcheon,
{\it An ergodic IP polynomial Szemer\'{e}di theorem},
Mem. Amer. Math. Soc. {\bf 146} (2000), no. 695.

\bibitem{BL93}
F. Blanchard and Y. Lacroix,
{\it Zero-entropy factors of topological flows},
Proc. Amer. Math. Soc. {\bf 119} (1993), no. 3, 985--992.


\bibitem{DDMSY13}
P. Dong, S. Donoso, A. Maass, S. Shao and X. Ye,
{\it Infinite-step nilsystems, independence and complexity},
Ergodic Theory Dynam. Systems {\bf 33} (2013), no. 1, 118--143.

\bibitem{FGJO21}
G. Fuhrmann, E. Glasner, T. J\"{a}ger and C. Oertel,
{\it Irregular model sets and tame dynamics},
Trans. Amer. Math. Soc. {\bf 374} (2021), no. 5, 3703--3734.


\bibitem{FH}
H. Furstenberg,
{\it Ergodic behavior of diagonal measures and a theorem of Szeme\'{r}edi on arithmetic progressions},
J. Anal. Math. {\bf 31} (1977), 204--256.


\bibitem{FH81}
H. Furstenberg,
{\it Recurrence in ergodic theory and combinatorial number theory},
M. B. Porter Lectures. Princeton University Press, Princeton, N. J. 1981.

\bibitem{GE94}
E. Glasner,
{\it Topological ergodic decompositions and applications to products of powers of a minimal transformation},
J. Anal. Math. {\bf 64} (1994), 241--262.


\bibitem{EG07}
E. Glasner,
{\it The structure of tame minimal dynamical systems},
Ergodic Theory Dynam. Systems {\bf 27} (2007), no. 6, 1819--1837.


\bibitem{GHSWY20}
E. Glasner, W. Huang, S. Shao, B. Weiss and X. Ye,
{\it Topological characteristic factors and nilsystems},
to appear in J. Eur. Math. Soc..


\bibitem{GW95}
E. Glasner and B. Weiss,
{\it Quasi-factors of zero-entropy systems},
J. Amer. Math. Soc. {\bf 8} (1995), no. 3, 665--686.

\bibitem{GY09}
E. Glasner and X. Ye,
{\it Local entropy theory},
Ergodic Theory Dynam. Systems {\bf 29} (2009), no. 2, 321--356.


\bibitem{HK05}
B. Host and B. Kra,
{\it Nonconventional ergodic averages and nilmanifolds},
Ann. of Math. (2) {\bf 161} (2005), no. 1, 397--488.


\bibitem{HK18}
B. Host and B. Kra,
{\it Nilpotent Structures in Ergodic Theory},
Mathematical surveys and monographs {\bf 236},
Providence, Rhode Island: American Mathematical Society, 2018.

\bibitem{HKM10}
B. Host, B. Kra and A. Maass,
{\it Nilsequences and a structure theory for topological dynamical systems},
Adv. Math. {\bf 224} (2010), no. 1, 103--129.


\bibitem{WH06}
W. Huang,
{\it Tame systems and scrambled pairs under an Abelian group action},
Ergodic Theory Dynam. Systems {\bf 26} (2006), no. 5, 1549--1567.

\bibitem{HLY}
W. Huang, H. Li and X. Ye,
{\it Localization and dynamical Ramsey property}, preprint.


\bibitem{HLSY03}
W. Huang, S. Li, S. Shao and X. Ye,
{\it Null systems and sequence entropy pairs},
Ergodic Theory Dynam. Systems {\bf 23}, (2003), no. 5, 1505--1523.


\bibitem{HLSY21}
W. Huang, Z. Lian, S. Shao and X. Ye,
{\it Minimal systems with finitely many ergodic measures},
J. Funct. Anal. {\bf 280} (2021), no. 12, Paper No. 109000, 42 pp.

\bibitem{HSY16}
W. Huang, S. Shao and X. Ye,
{\it Nil Bohr-sets and almost automorphy of higher order},
Mem. Amer. Math. Soc. {\bf 241} (2016) no. 1143.


\bibitem{HSY17}
W. Huang, S. Shao and X. Ye,
{\it Strictly ergodic models under face and parallelepiped group actions},
Commun. Math. Stat. {\bf 5} (2017), no. 1, 93--122.

\bibitem{HY06}
W. Huang and X. Ye,
{\it A local variational relation and applications},
Israel J. of Math. {\bf 151} (2006), 237--280.

\bibitem{KAS95}
A.S. Kechris,
{\it Classical Descriptive Set Theory, Graduate Texts in Mathematics},
vol. 156, Springer-Verlag, New York, 1995.


\bibitem{KL07}
D. Kerr and H. Li,
{\it Independence in topological and $C^*$-dynamics},
Math. Ann. {\bf 338} (2007), no. 4, 869--926.


\bibitem{KL16}
D. Kerr and H. Li,
{\it Ergodic theory. Independence and dichotomies},
Springer Monographs in Mathematics. Springer, Cham, 2016.

\bibitem{KK66}
K. Kuratowski,
{\it Topology}, Vol. I, Acad. Press, New York, N.Y., 1966.

\bibitem{KK68}
 K. Kuratowski,
{\it Topology}, Vol. II, Acad. Press, New York, N.Y., 1968.

\bibitem{MS07}
A. Maass and S. Shao,
{\it Structure of bounded topological sequence entropy minimal systems},
J. Lond. Math. Soc. (2) {\bf 76} (2007), no. 3, 70--718.

\bibitem{SY12}
S. Shao and X. Ye,
{\it Regionally proximal relation of order d is an equivalence one for minimal systems and a combinatorial consequence},
Adv. Math. {\bf 231} (2012), no. 3--4, 1786--1817.

\bibitem{Sz75}
E. Szemer\'{e}di,
{\it On sets of integers containing no k elements in arithmetic progression},
Acta Arith. {\bf 27} (1975), 199--245.


\bibitem{JDV}
J. de Vries,
{\it Elements of Topological Dynamics},
Kluwer Academic Publishers (993), Dordrecht.

\bibitem{TZ07}
T. Ziegler,
{\it Universal characteristic factors and Furstenberg averages},
J. Amer. Math. Soc. {\bf 20} (2007), no.1, 53--97.

\end{thebibliography}
\end{document}